\documentclass[12pt]{amsart}
\usepackage{epsfig,amsmath,amsfonts,latexsym,accents}
\usepackage[usenames,dvipsnames]{color}
\usepackage{palatino}
\usepackage[hypcap=false]{caption}
\usepackage{caption, subcaption, pinlabel, float}
\usepackage{color}
\usepackage[normalem]{ulem}

\usepackage{xcolor}
\usepackage{comment}

\definecolor{darkgreen}{rgb}{0, 0.5, 0}
\definecolor{darkpurple}{rgb}{.5,0,.5}
\definecolor{darkbrown}{rgb}{.62,.35,.17}
\usepackage[hyperfootnotes=false]{hyperref}
\hypersetup{
  colorlinks,
  citecolor=Purple,
  linkcolor=Black,
  urlcolor=Green}

\theoremstyle{plain}
\newtheorem{dummy}{anything}[section]
\newtheorem{theorem}[dummy]{Theorem}

\newtheorem{lemma}[dummy]{Lemma}

\newtheorem{question}[dummy]{Question}

\newtheorem{proposition}[dummy]{Proposition}

\newtheorem{corollary}[dummy]{Corollary}

\theoremstyle{definition}
\newtheorem{definition}[dummy]{Definition}

\newtheorem{example}[dummy]{Example}

\newtheorem{remark}[dummy]{Remark}

\theoremstyle{remark}

\newcommand{\del}{\partial}

\newcommand{\Z}{\mathbb{Z}}
\newcommand{\R}{\mathbb{R}}
\newcommand{\C}{\mathbb{C}}

\newcommand{\K}{\rm{K}}

\newcommand\ttimes{\mathbin{%
    \stackrel{\sim}{\smash{\times}\rule{0pt}{0.7ex}}%
    }}

\newcommand{\bP}{\mathbb P}

\def\bdy{\partial}

\def\a{\alpha}
\def\b{\beta}

\def\g{\gamma}

\begin{document}

\title {Lefschetz fibrations on nonorientable 4-manifolds}

\author{Maggie Miller}

\author{Burak Ozbagci}

\address{Department of Mathematics, Massachusetts Institute of Technology, Cambridge,
MA 02142, USA}

\email{maggiehm@mit.edu}

\address{Department of Mathematics, Ko\c{c} University, Istanbul,
Turkey}
\email{bozbagci@ku.edu.tr}

\subjclass[2000]{}
\thanks{}


\begin{abstract}
Let $W$ be a nonorientable $4$-dimensional handlebody without $3$- and $4$-handles. We show that $W$ admits a Lefschetz fibration over the $2$-disk, whose regular fiber is a nonorientable surface with nonempty boundary.  This is an analogue of a result of Harer obtained in the orientable case.  As a corollary,  we obtain a $4$-dimensional proof of the fact that every nonorientable closed  $3$-manifold admits an open book decomposition, which was first proved by Berstein and Edmonds using branched coverings.   Moreover, the monodromy of the open book we obtain for a given $3$-manifold belongs to the twist subgroup of the mapping class group of the page. In particular, we construct an explicit minimal open book for the connected sum of arbitrarily many copies of the product of the circle with the real projective plane.

We also obtain a relative trisection diagram for $W$, based on the nonorientable Lefschetz fibration we construct, similar to the orientable case first studied by Castro. As a corollary, we get trisection diagrams for some closed $4$-manifolds, e.g. the product of the $2$-sphere with the real projective plane,  by doubling $W$.   Moreover, if $X$ is a  closed nonorientable $4$-manifold  which admits a  Lefschetz fibration over the $2$-sphere, equipped with a section of square $\pm 1$,  then we construct a trisection diagram of $X$,  which is determined by the vanishing cycles of the Lefschetz fibration. Finally, we include some simple observations about low-genus Lefschetz fibrations on closed nonorientable $4$-manifolds.

 \end{abstract}

\maketitle

\section{Introduction}\label{sec: intro}  Since Lefschetz fibrations were put into the foreground at the turn of this century by the seminal works of Donaldson and Gompf,  various flavors (e.g. achiral, broken, bordered, symplectic, holomorphic) of Lefschetz fibrations have been fruitfully utilized in order to study some aspects of smooth {\em orientable} manifolds, yielding especially interesting results in dimension $4$.

In this paper, we initiate the study of Lefschetz fibrations on compact  {\em nonorientable} $4$-manifolds. We have opted to take the base of any of these fibrations as a compact orientable $2$-manifold so that a regular fiber is necessarily a nonorientable surface. We say that a Lefschetz fibration is of genus $g$ if the (nonorientable) fiber is of genus $g$ as a surface. Here, the genus of a closed nonorientable surface is defined to be the (positive) number of crosscaps attached to $S^2$ to obtain the surface at hand. 

When the nonorientable $4$-manifold at hand is closed, one can explore the interaction between the given Lefschetz fibration and the mapping class group of the  nonorientable regular fiber of the fibration. A notable distinction from the well-studied orientable case is that Dehn twists  generate the  index two twist subgroup  of  the mapping class group of a nonorientable surface  \cite{l}.

On the other hand, a Lefschetz fibration over $D^2$ on a compact nonorientable $4$-manifold with nonempty boundary naturally induces an open book on the boundary $3$-manifold, which is also nonorientable.  For example,  a genus one Lefschetz fibration $D^2 \times \R\bP^2 \to D^2$ we construct  in Example~\ref{ex: prod} led us to the discovery of a genus one open book for $S^1 \times \R\bP^2$. As a consequence, for each $n \geq 2$, we construct in Example~\ref{ex: connect} an explicit genus one open book for $\#_n S^1 \times \R\bP^2$ by successively performing Murasugi sums.

Our first result in this direction is the nonorientable analogue of a result proved by Harer  \cite{h}  for the orientable case.

\begin{theorem} \label{thm: lef} Let $W$ be a nonorientable $4$-dimensional handlebody without $3$- and $4$-handles. Then $W$ admits an explicit Lefschetz fibration over $D^2$, whose regular fiber is a nonorientable surface with nonempty boundary.
\end{theorem}

As an immediate corollary to Harer's result, one obtains a proof of a classical theorem of Alexander \cite{a} that says: {\em Every orientable closed $3$-manifold admits an open book decomposition.} In the same vein, we obtain  Corollary~\ref{cor: openb}, which was first proved by Berstein and Edmonds \cite[Theorem 9.8]{be} using branched coverings, but we would like to emphasize that the additional information about the monodromy, namely that it  can be expressed as a product of Dehn twists,  does not follow from their construction.

\begin{corollary}\label{cor: openb}
Every nonorientable closed  $3$-manifold admits an open book decomposition, whose monodromy can be expressed as a product of Dehn twists.
\end{corollary}

As a matter of fact,   Ghanwat, Pandit and Selvakumar \cite{gps} recently showed that every closed nonorientable $3$-manifold admits a {\em genus one} open book (meaning that its page is the real projective plane $\R\bP^2$ with holes), whose monodromy is a product of Dehn twists. Moreover, they described an algorithm to compute the monodromy of the open book they construct for a given $3$-manifold. Their procedure involves writing a $3$-manifold as Dehn surgery on a link which is braided in a specific way about the binding of a genus one open book for $S^2\ttimes S^1$.

In Example~\ref{ex: genusone} below,  we obtain an explicit genus one open book for  $S^1 \times \R\bP^2$ based on our Theorem~\ref{thm: lef}. Reading \cite{gps} literally, their procedure applied to $S^1 \times \R\bP^2$ yields an open book with three binding components, but by starting with the braid in Figure \ref{fig:ghanwat}, we see that their procedure on braids yields the same open book depicted in Figure~\ref{fig:introfig}, which is obtained via our Theorem~\ref{thm: lef}. We should also point out that Klassen \cite{k} constructed an explicit genus two open book for $S^1 \times \R\bP^2$ (described somewhat less explicitly in \cite{be}), whose monodromy is the $Y$-homeomorphism of Lickorish. The interested reader may consult \cite{o} for further details. Berstein and Edmonds \cite{be} use this open book to construct an open book for an arbitrary nonorientable closed $3$-manifold $M$ by writing $M$ as an $n$-fold branched cover over $S^1\times\R\bP^2$ and lifting this open book to $M$. When $n$ is even, the resulting open book has monodromy that can be written as a product of Dehn twists; when $n$ is odd then not.

\begin{remark} A simple closed curve $\g$ in a nonorientable surface is called two-sided if a regular
neighborhood of $\g$  is an annulus. A Dehn twist about $\g$ is defined as usual,  but
to be able to distinguish a  right-handed twist from a left-handed twist, one must also specify  an orientation of a regular
neighborhood of $\g$, for which there is no canonical choice.
\end{remark}

\begin{example} \label{ex: genusone} There is a genus one open book for $S^1 \times \R\bP^2$ whose page is  $\R\bP^2$ with two holes and whose monodromy is the product of Dehn twists about the boundary parallel  curves $\gamma_1$ and $\gamma_2$ depicted in Figure~\ref{fig:introfig}.

\begin{figure}[H]{\centering
\vspace{.05in}
\labellist
  \pinlabel {\footnotesize\textcolor{darkgreen}{$\gamma_1$}} at 41 80
  \pinlabel {\footnotesize\textcolor{darkgreen}{$\gamma_2$}} at 57 30
\endlabellist
\includegraphics[width=1.5in]{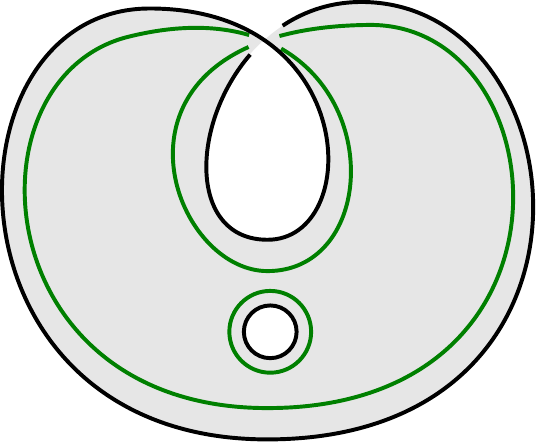}

\caption{A page of an open book decomposition of $S^1\times\R\bP^2$. The monodromy is a product of Dehn twists about $\gamma_1$ and $\gamma_2$.}
\label{fig:introfig}
}
\end{figure}

This follows immediately from the genus one Lefschetz fibration $D^2 \times \R\bP^2 \to D^2$ we construct in Example~\ref{ex: prod}, since we have 
$\bdy(D^2 \times  \R\bP^2)= S^1 \times \R\bP^2$. Alternatively, as we discuss with more details in  Remark~\ref{rem: pi1}, one can directly check that the total space of the abstract open book whose page is  $\R\bP^2$ with two holes and whose monodromy is the product of the two boundary Dehn twists (with arbitrary ``handedness" for a fixed choice of orientations of the annuli neighborhoods of  $\gamma_1$ and $\gamma_2$, see Figure \ref{fig:dehnsymmetry})  is diffeomorphic to $S^1 \times \R\bP^2$. \end{example}

\begin{figure}{\centering
\vspace{.05in}
\labellist
  \pinlabel {\footnotesize\textcolor{darkgreen}{$1$}} at 95 65
  \pinlabel {\footnotesize\textcolor{darkgreen}{$-1$}} at 160 115
  \pinlabel {\footnotesize\textcolor{darkgreen}{$1$}} at 185 65
  \pinlabel {\footnotesize\textcolor{darkgreen}{$-1$}} at 250 115
\endlabellist
\includegraphics[width=4.5in]{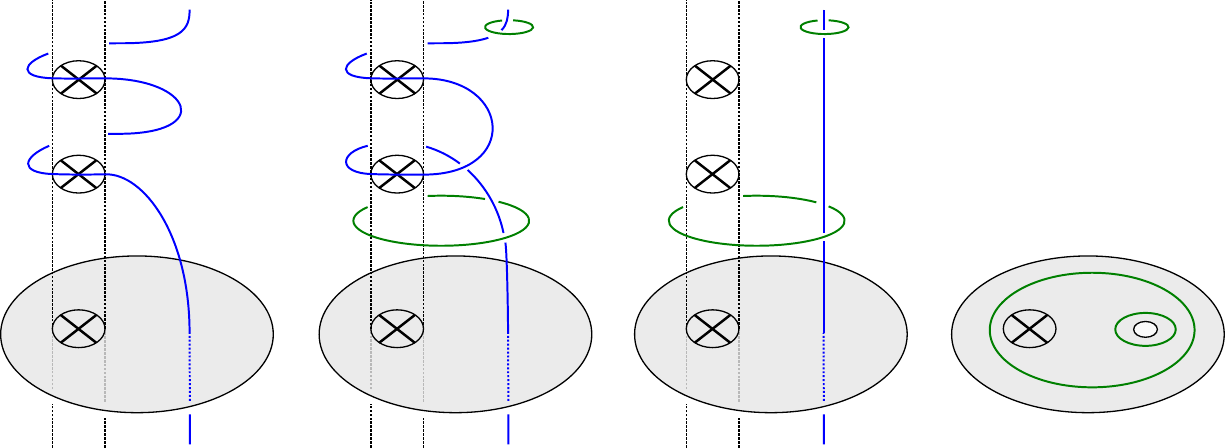}

\caption{The procedure of \cite{gps} to obtain a genus one open book on $S^1\times\R\bP^2$: On the left, we draw a M\"{o}bius band $M$ representing the page of an open book on $S^2\ttimes S^1$. We write $S^1\times \R\bP^2$ as Dehn surgery on a link (in this case a knot) braided about $M$ as in \cite{gps}. We blow up and isotope (second and third image) to make the surgery link into a 0-framed vertical component and two $\pm 1$ unknotted components parallel to the M\"{o}bius band. Finally on the right, we puncture the M\"{o}bius band at the vertical component to obtain the page of a genus one open book on $S^1\ttimes\R\bP^2$. The monodromy consists of Dehn twists about the two green curves (projections of the unknotted surgery curves) with opposite handedness (due to the signs of the two surgery curves being opposite.) Unlabeled arcs are all 0-framed.}
\label{fig:ghanwat}
}
\end{figure}
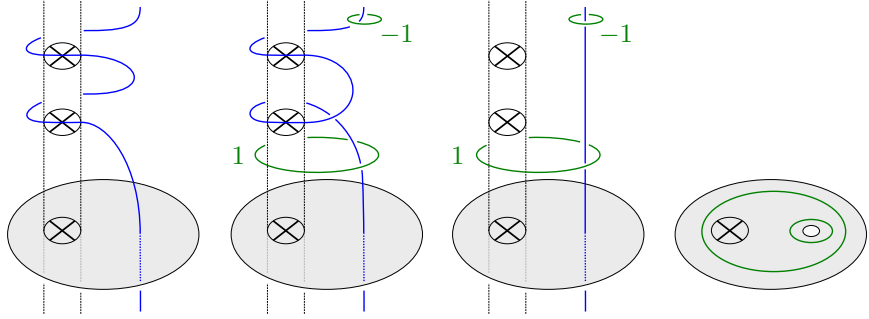

\begin{figure}[H]{\centering
\vspace{.05in}
\labellist
  \pinlabel {\footnotesize\textcolor{darkgreen}{$\gamma_2$}} at 68 116
  \pinlabel {\footnotesize\textcolor{darkgreen}{$\gamma_1$}} at 45 50
\endlabellist
\includegraphics[width=1.5in]{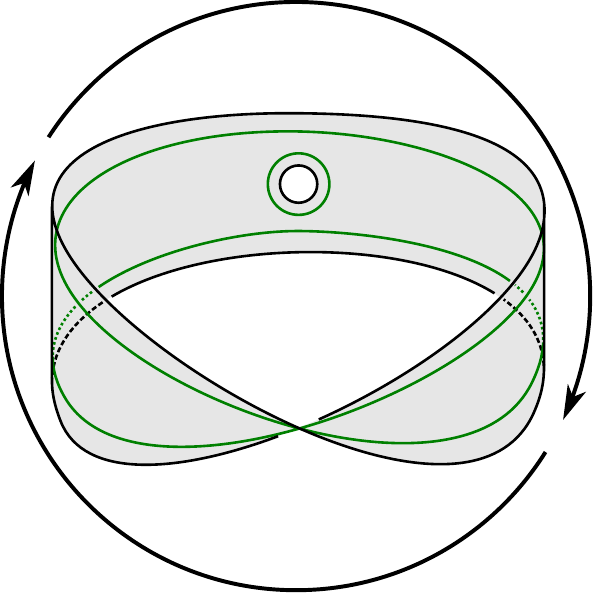}

\caption{The ''handedness'' of each Dehn twist in the open book decomposition of $S^1\times\R\bP^2$ depicted in  Figure \ref{fig:introfig} does not matter, i.e. any choice gives an equivalent open book. Here we give a diffeomorphism of the page that fixes each of $\gamma_1$ and $\gamma_2$ setwise but is orientation-reversing on $\gamma_1$ and orientation-preserving on $\gamma_2$. (The diffeomorphism is a 180$^\circ$ rotation in the pictured ambient 3D space as indicated by the arrows.) This gives an equivalence between two open books with monodromy a Dehn twist on $\gamma_1$ and $\gamma_2$ with the opposite choice of twist on $\gamma_1$ and the same choice of twist on $\gamma_2$.}
\label{fig:dehnsymmetry}
}
\end{figure}

\begin{example} \label{ex: connect} Relying on the nonorientable version (cf. \cite[Proposition 6.1]{o})  of Stallings' Theorem \cite{s} and performing  a Murasugi sum of two copies of the genus one open book for  $S^1\times\R\bP^2$ depicted in Figure~\ref{fig:introfig},  we get a genus one open book for $\#_2 S^1 \times \R\bP^2$,  whose page is  $\R\bP^2$ with four holes and whose monodromy is the product of Dehn twists (with arbitrary handedness) about the curves $\gamma_1, \gamma_2, \gamma_3$ and $\gamma_4$ depicted in Figure~\ref{fig:genus1}. 

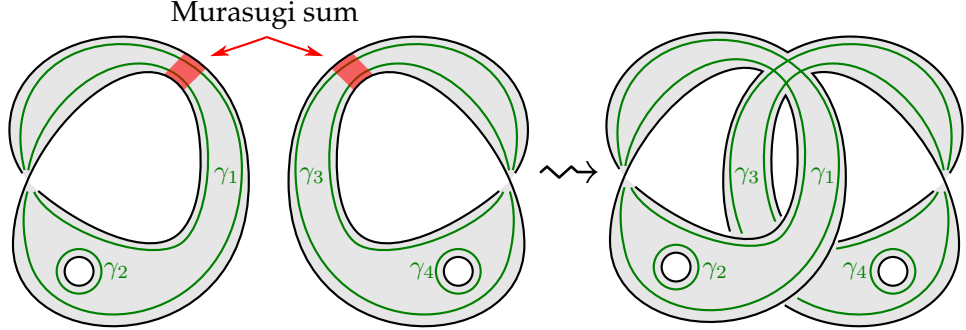
\begin{figure}[H]{\centering
\vspace{.15in}
\labellist
\pinlabel {\small{Murasugi sum}} at 135 165
\pinlabel {\huge{$\rightsquigarrow$}} at 295 80
  \pinlabel {\footnotesize\textcolor{darkgreen}{$\gamma_1$}} at 115 80
  \pinlabel {\footnotesize\textcolor{darkgreen}{$\gamma_2$}} at 57 30
   \pinlabel {\footnotesize\textcolor{darkgreen}{$\gamma_3$}} at 160 80
    \pinlabel {\footnotesize\textcolor{darkgreen}{$\gamma_4$}} at 218 30
      \pinlabel {\footnotesize\textcolor{darkgreen}{$\gamma_1$}} at 429 80
  \pinlabel {\footnotesize\textcolor{darkgreen}{$\gamma_2$}} at 372 30
   \pinlabel {\footnotesize\textcolor{darkgreen}{$\gamma_3$}} at 389 80
    \pinlabel {\footnotesize\textcolor{darkgreen}{$\gamma_4$}} at 446 30
\endlabellist
\includegraphics[width=5in]{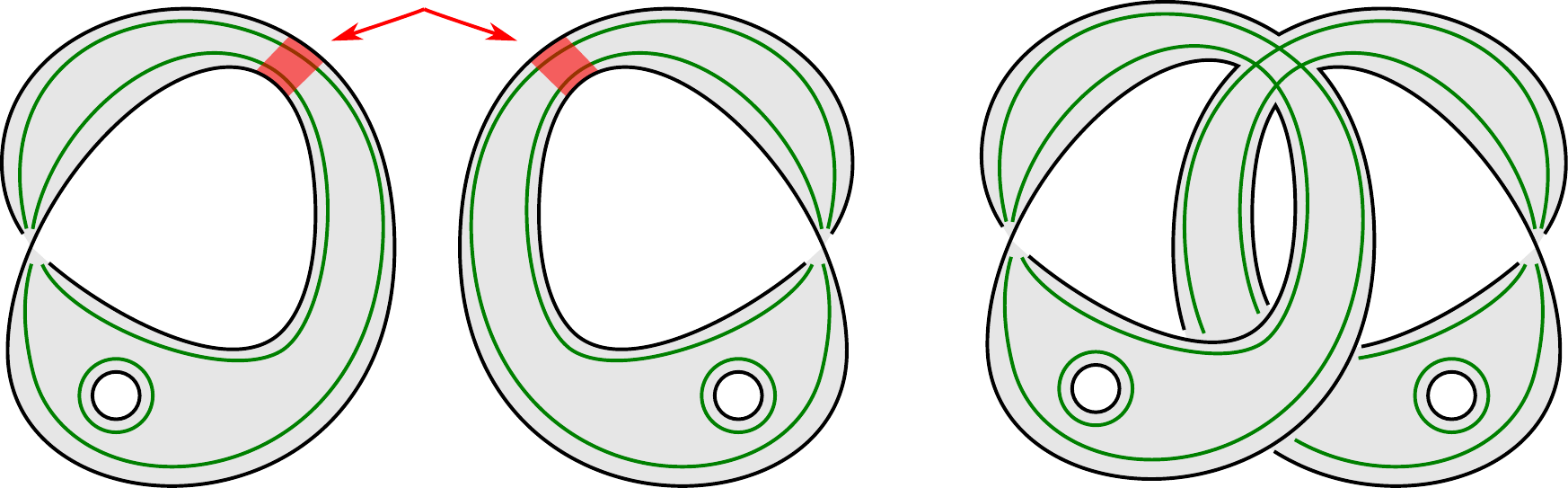}

\caption{We construct a genus one open book for $\#_2 S^1\times\R\bP^2$ via Murasugi sum by identifying the colored rectangles as indicated on the left.  The monodromy of the resulting open book on the right is given by a product of four Dehn twists: one each about the curves  $\gamma_1$, $\gamma_2$, $\gamma_3$, $\gamma_4$, respectively, in any order. Notice that the order of Dehn twists about $\gamma_2$ and $\gamma_4$ does not matter since they commute with the others, and the order of Dehn twists about $\gamma_1$ and $\gamma_3$ can be interchanged by an overall conjugation of the monodromy.  }
\label{fig:genus1}
}
\end{figure}

Similarly, for 
$n \geq 3$, we may Murasugi sum $n$ copies of the genus one open book for $S^1\times\R\bP^2$ (iteratively plumbing $n-1$ copies to one base copy along the rectangles indicated in Figure \ref{fig:genus1more}) to obtain a genus one open book for $\#_n S^1\times\R\bP^2$ whose page is $\R\bP^2$ with $2n$ holes. For any 
$n \geq 2$, the genus one open book we construct for $\#_n S^1\times\R\bP^2$ is in fact induced by a Lefschetz fibration $\natural_n D^2 \times \R\bP^2 \to D^2$, where $\natural$ denotes the boundary connected sum, which follows from a nonorientable version of a result in  \cite{koo}. \end{example}

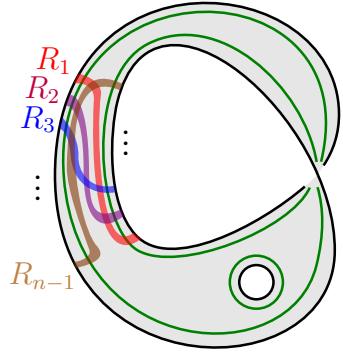
\begin{figure}{\centering
\labellist
  \pinlabel {\textcolor{red}{$R_1$}} at 0 128
  \pinlabel {\textcolor{purple}{$R_2$}} at -5 115
   \pinlabel {\textcolor{blue}{$R_3$}} at -7 102
   \pinlabel {$\vdots$} at -7 75
     \pinlabel {$\vdots$} at 32 95
    \pinlabel {\textcolor{brown}{$R_{n-1}$}} at -5 32
\endlabellist
\includegraphics[width=1.5in]{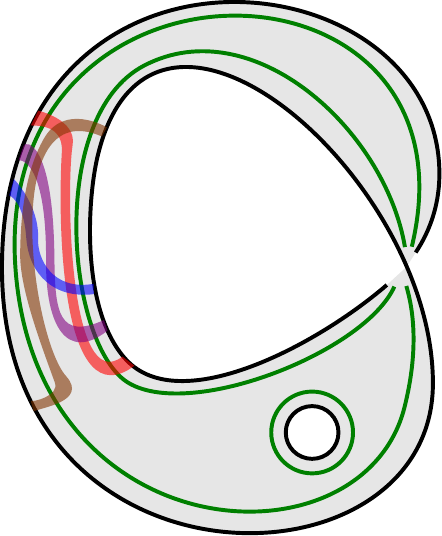}

\caption{{Murasugi summing $n-1$ copies of the genus one open book for $S^1\times\R\bP^2$ (as in Figure \ref{fig:introfig}) to the pictured copy along rectangles $R_1,\ldots, R_{n-1}$ in order yields a genus one open book for $\#_n S^1\times\R\bP^2$. (The plumbing rectangles on the unpictured copies are as in Figure \ref{fig:genus1}.) {The resulting page is a M\"obius band with $2n-1$ holes. } }
\label{fig:genus1more}}}
\end{figure}

Proposition~\ref{prop: minimal} implies that the genus one open books we constructed in Example \ref{ex: connect} (and also Example \ref{ex: genusone}) are all ``minimal" in the sense that  the page is of genus one (minimum possible) and the number of binding components is minimal for a genus one open book.

\begin{proposition}\label{prop: minimal}
For any $n\ge 1$, a genus one open book for $\#_n S^1\times\R\bP^2$ must have at least $2n$ binding components.
\end{proposition}

As a byproduct of the construction in our proof of Theorem~\ref{thm: lef}, we obtain the following result, which might be of interest on its own.

\begin{proposition}\label{prop: embed} If $L$ is a framed link in $\#_p S^2 \ttimes S^1$, for some nonnegative integer $p$,  then there is an explicit open book decomposition for $\#_p S^2 \ttimes S^1$ such that each component of $L$ is embedded in a distinct page with framing $\pm1$ relative to the page framing. \end{proposition}

\begin{remark}\label{frameremark}
Knot framings are slightly more subtle in nonorientable manifolds, but there is no subtlety in relative framings. Let $C$ be a simple closed curve in a surface $\Sigma$ that is a page of an open book in a 3-manifold $M$. In order to consider framings of $C$, $C$ must have trivial normal bundle, so we restrict to the case that $C$ has an annular neighborhood in $\Sigma$. Then a tangent and normal vector to $\Sigma$ at $C$ induce a framing $\xi$ of the normal bundle of $C$ in $M$. We say that the framing $\eta$ of $C$ differs from $\xi$ by $\pm 1$ if $\eta$ can be isotoped to agree with $\xi$ away from a small interval in $C$, where they differ by a single twist of appropriate sign. If $M$ is nonorientable, then the sign does not matter: if we can isotope $\xi$ to agree with $\eta$ except for one positive twist, then we can further isotope $\xi$ to agree with $\eta$ except for one negative twist.
\end{remark}

Next, we  turn our attention to the interaction between Lefschetz fibrations and trisections. In his PhD thesis \cite{c}, Castro constructed a relative trisection of any compact orientable $4$-manifold which admits an achiral Lefschetz fibration over $D^2$. Here we observe that his result extends to the nonorientable case to obtain another corollary of Theorem~\ref{thm: lef}.

\begin{corollary}\label{cor: reltri}  Let $W$ be a nonorientable $4$-dimensional handlebody without $3$- and $4$-handles. Then, based on any Lefschetz fibration $W \to D^2$,  there is an explicit algorithm to obtain a relative trisection diagram of $W$ so that the open book on $\bdy W$ induced by the Lefschetz fibration coincides with that induced by the trisection.
\end{corollary}

As an application, in Remark~\ref{rem:yhomeo} below, we show that there exists a compact nonorientable $4$-manifold without $3$- and $4$-handles, which admits two distinct relative trisections that are not stably equivalent  via relative trisection boundary stabilizations as in \cite{c}.

Let $X$ be a closed, nonorientable $4$-manifold and let $X_{(2)}$ denote the union of the $0$-, $1$-, and $2$-handles in a given handle decomposition of $X$. Then, as shown in the recent work of the first author and Naylor \cite[Corollary 3.13]{mn},   $X$ is determined up to diffeomorphism by $X_{(2)}$, which follows from a  nonorientable analogue of a classical result of  Laudenbach and Po\'{e}naru \cite{lp}. Therefore, the relative trisection diagram of $X_{(2)}$ described in Corollary~\ref{cor: reltri} determines the diffeomorphism type of the {\em closed} $4$-manifold $X$.

Nevertheless, if desired, one can construct trisection diagrams for some closed nonorientable $4$-manifolds based on Corollary~\ref{cor: reltri} as follows. Let $X$ be the double of $W$, i.e., $X = W \cup_\partial W$. Then the identical relative trisection diagrams on the two copies of $W$ can be glued together to obtain a trisection diagram of the closed nonorientable $4$-manifold $X$. We illustrate this method in Example~\ref{ex: s2rp2}.

Moreover, if  $X$ is a  closed, nonorientable, connected $4$-manifold  which admits a Lefschetz fibration over $S^2$  equipped with a section of square $\pm 1$, then we show that $X = W \cup_{\bdy} V$, where both $W$ and $V$ admit Lefschetz fibrations over $D^2$ so that the induced open books coincide on their common boundary $\bdy W = \bdy V$.  Using Corollary~\ref{cor: reltri}, we obtain relative trisection diagrams  of $W$ and $V$, respectively and glue these diagrams to get a trisection diagram of $X$.  We advice the reader to turn to \cite{cgpc, co}  for a detailed discussion on gluing relative trisection diagrams in the orientable case and to \cite{mn} in the nonorientable case. The discussion in this paragraph is summarized as Theorem~\ref{thm: closedLEF} below, which is a nonorientable analogue of Theorem 3.7 in the article \cite{co} of Castro and the second author.

\begin{theorem} \label{thm: closedLEF} Suppose that $X$ is a  closed, nonorientable, connected $4$-manifold  which admits a  Lefschetz fibration over $S^2$  equipped with a section of square $\pm 1$. Then, an explicit  trisection of $X$ can be described  by a  corresponding trisection diagram,  which is determined by the vanishing cycles of the Lefschetz fibration.    \end{theorem}

Finally, we include a simple observation which is a nonorientable analogue of the following fact (see e.g. \cite[Proposition 8.1.7]{gs}): ``{\em A relatively minimal genus zero Lefschetz fibration on a closed oriented $4$-manifold over any closed oriented surface is an $S^2$-bundle.}" (See Section~\ref{sec: harer}, for the relevant definitions.)

\begin{proposition} \label{prop: genoneLEF} If $X \to B$ is a relatively minimal genus one Lefschetz fibration on a closed nonorientable $4$-manifold $X$ over a closed orientable surface $B$,  then $X$ is an $\R\bP^2$-bundle over $B$.  
\end{proposition} 

In other words,  there are no singular fibers in a relatively minimal genus one Lefschetz fibration  $X \to B$, where $X$ is  any closed nonorientable $4$-manifold and $B$ is an arbitrary orientable surface.  Note that there are two diffeomorphism classes of nonorientable total spaces of $\R\bP^2$-bundles over an orientable surface $B$, classified by the second Stiefel-Whitney class of the bundle.

On the other hand, in Section~\ref{sec: closed}, for any $g \geq 2$, we construct examples of  relatively minimal (nonorientable) genus $g$ Lefschetz fibrations over $B$ with arbitrary number of singular fibers.

\section{A nonorientable analogue of Harer's result}\label{sec: harer}
The standard definition of a Lefschetz fibration on orientable $4$-manifolds (cf.  \cite[Section 8]{gs}) can be easily adapted to nonorientable $4$-manifolds.

\begin{definition} \label{def: lefs} Let $X$ be a compact, connected, {\em nonorientable} $4$-manifold and let $B$ be a compact, connected,   {\em orientable} $2$-manifold.   A Lefschetz fibration is a surjective map $\pi : X \to B$ with finitely many critical points in the interior of $X$ such that
around each critical point, $\pi $ conforms to the model $ \pi (z_1, z_2) = z_1z_2$ in local complex coordinates. 
\end{definition}

The manifolds $X$ and $B$ in Definition~\ref{def: lefs} are allowed to have nonempty boundaries. In that case, if the regular fiber is a closed surface, then $\bdy X = \pi^{-1}(\bdy B)$ and if the regular fiber is a surface with nonempty boundary, then  $\bdy X$ acquires a natural open book whose page is the regular fiber.

A Lefschetz fibration is called {\em relatively minimal} if no fiber contains an exceptional sphere, i.e., a sphere with self-intersection $\pm 1$. Any Lefschetz fibration can be blown down to obtain a relatively minimal one. 

\begin{remark} \label{rem: comp} Since $X$ is {\em nonorientable}, we can not impose the local complex coordinates around a critical point to be compatible with the ``orientation" of $X$. Moreover,  a regular fiber $F$ is a (necessarily) nonorientable surface such that $\bdy F \neq \emptyset$ if and only if $\bdy X \neq \emptyset$.  We may assume that each
singular fiber contains only one critical point, which can be achieved after a small perturbation of $\pi$. \end{remark}

In this section, we will be interested only in the case when $B=D^2$. Suppose that $\pi: X \to D^2$ is a Lefschetz fibration. It follows that each fiber of $\pi$  is connected (cf. \cite[Proposition 8.1.9]{gs}). Let $F$ denote a regular fiber of the Lefschetz fibration  $\pi: X \to D^2$. The topology of $X$ can be described as follows. Start with the trivial fibration $D^2\times F $ and attach $2$-handles along {\em two-sided} simple closed curves---one for each critical point---embedded in distinct regular fibers along $\bdy (D^2\times F )$. The framing of each  $2$-handle is given by $\pm 1$ with respect to the surface framing. Notice that one can not distinguish $+$ or $-$ because, although  a regular neighborhood of the attaching curve on $F$ is orientable, there is no canonical orientation one can choose. Similar to the orientable Lefschetz fibrations, the local monodromy around a critical value is a Dehn twist around "the vanishing cycle" with the caveat that there is no canonical choice of right-handed versus left-handed in the nonorientable case.   We are now ready to give a proof of  Theorem~\ref{thm: lef}, adapting  the proof in Harer's Thesis \cite{h} to the nonorientable case. See also \cite{ef}.

\begin{proof} {\em of Theorem~\ref{thm: lef}.}  Let $W$ be a nonorientable $4$-dimensional handlebody without $3$- and $4$-handles. Then $W$ may be obtained by attaching  $2$-handles to $W_{(1)} \cong \natural_p D^3 \ttimes S^1$ for some positive integer $p$. In other words, $W$ consists of a single $0$-handle, $p$ nonorientable $1$-handles and some $2$-handles. Notice that $$W_{(1)} \cong \natural_p D^3 \ttimes S^1 \cong D^2\times \Sigma,$$ where $\Sigma$ is the nonorientable surface obtained by attaching $p$ nonorientable $1$-handles to $D^2$. We can think of $  D^2\times\Sigma$ as the trivial Lefschetz fibration over $D^2$. More importantly,  we can view $\Sigma$ as a disk-with-twisted bands (see left of Figure \ref{fig:fig11}) inside  $W_{(1)} \subset  W$, which is a page of the trivial open book for $\#_p S^2 \ttimes S^1$ induced by the trivial Lefschetz fibration  $W_{(1)} \cong D^2\times \Sigma  \to D^2$.

Next, we would like to embed each component of the link $L \subset \#_p S^2 \ttimes S^1$, consisting of the attaching curves of the $2$-handles,  in a distinct page of this open book. To this end, we take a generic projection of  $L$ onto $\Sigma$ with only transverse double points and then apply {\em Harer's trick} as follows. Take the connected sum of $\Sigma$ with a torus at each double point of the projection as shown in Figure \ref{fig:fig11} (middle).

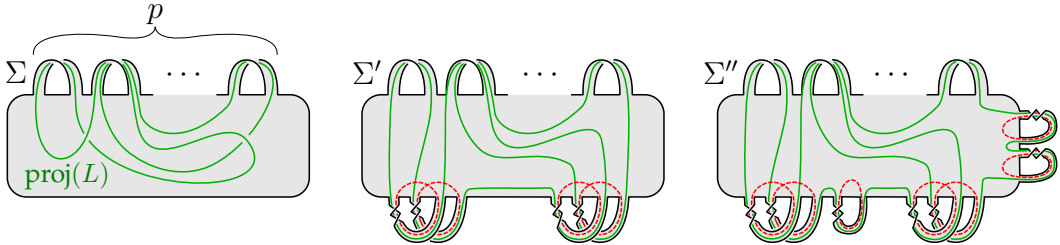
\begin{figure}[H]{\centering
\vspace{.05in}
\labellist
  \pinlabel {\footnotesize\textcolor{darkgreen}{proj$(L)$}} at 45 50
   \pinlabel {$\Sigma$} at 7 127
     \pinlabel {$\Sigma'$} at 267 127
        \pinlabel {$\Sigma''$} at 529 127
    \pinlabel {$\cdots$} at 133 125
      \pinlabel {$\cdots$} at 397 125
        \pinlabel {$\cdots$} at 657 125
    \pinlabel {$p$} at 110 170
\endlabellist
\includegraphics[width=5.5in]{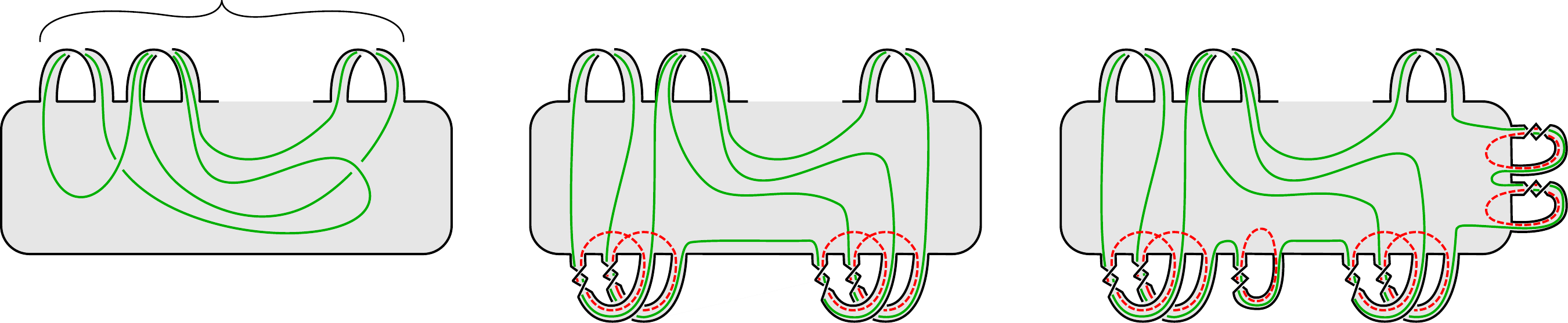}

\caption{Compare to Figures 2 and 3 of \cite{ef}. {\bf{Left:}} the link $L$ projected onto $\Sigma$, as in Theorem \ref{thm: lef}. {\bf{Middle:}} We perform Harer's trick at the double points of the projection of $L$. We indicate new vanishing cycles in red (dashed). {\bf{Right:}} We further stabilize $\Sigma$ to change the induced framing on $L$. Note that at each stabilization, we add a whole twist to $L$ on the stabilized surface.}
\label{fig:fig11}
}
\end{figure}

This modification corresponds to stabilizing the Lefschetz fibration (and thus the open book on the boundary) twice, keeping the diffeomorphism type of $W$ fixed. Once all the double points of the projection of $L$ is taken care of, we end up with a (nontrivial) Lefschetz fibration $W \to D^2$ such that each component $L_i \subset L$ is embedded in a distinct page $\Sigma'$ (obtained by stabilizing $\Sigma$ by several times)  of the induced open book on $\bdy W$.

Since $L$ consists of the attaching circles of the $2$-handles, each component of $L$ has trivial normal bundle in $W$. Since $\Sigma'$ is the page of an open book, the normal bundle of $\Sigma'$ in $W$ is also trivial. Therefore, each component of $L$ has an annular neighborhood in $\Sigma'$. We adjust the framings of the components of $L$ so that the framing of each component $L_i \subset L$ is $\pm  1$ with respect to the page framing (see Remark \ref{frameremark}). This is easily accomplished by stabilizing the page  $\Sigma'$ further as shown in the right of Figure \ref{fig:fig11} to obtain $\Sigma''$ which carries additional vanishing cycles corresponding to each stabilization. To summarize, we have now a Lefschetz fibration $W \to D^2$ so that each $L_i \subset L$ is embedded in a distinct fiber in the boundary, with framing $\pm  1$ with respect to the fiber framing.  Once the $2$-handles are attached along the link $L$, we get the desired Lefschetz fibration on $W$, whose fiber is $\Sigma''$ and each $L_i$ being a vanishing cycle in addition to the cycles inserted for each stabilization in the construction. 

The ``handedness" of the vanishing cycle corresponding to $L_i$ is determined by the framing of $L_i$: if we arrange the normal bundle of $L_i$ to have 2-handle framing agreeing with that induced by $\Sigma''$ except for one extra positive twist (with respect to some local orientation of $\Sigma''$) of the 2-handle attachment in a small arc, then near that same arc the corresponding Dehn twist in the monodromy induced on $\Sigma''$ is right-handed (with respect to the same local orientation of $\Sigma''$). \end{proof}

\begin{figure}{\centering
\subcaptionbox{A handlebody diagram of $D^2\ttimes\R\bP^2$.
\label{fig:twdiskbunda}}{
\labellist
  \pinlabel {$1$} at 61 18
\endlabellist
\includegraphics[width=1.5in]{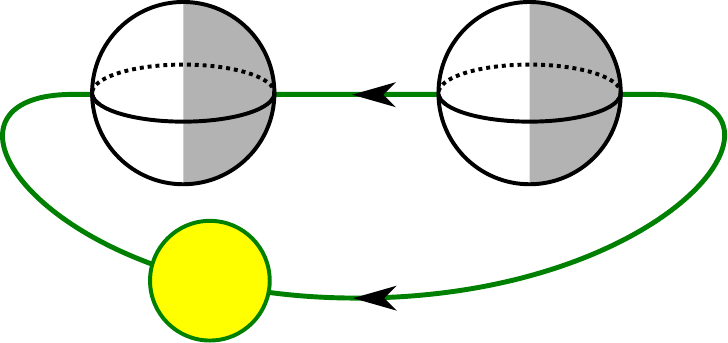}
}
\hspace{1cm}
\subcaptionbox{$D^2\ttimes\R\bP^2$ admits a genus one Lefschetz fibration over $D^2$ with a single vanishing cycle $\gamma$.
\label{fig:twdiskbundb}}{
\labellist
  \pinlabel {\textcolor{darkgreen}{$\gamma$}} at 61 22
\endlabellist
\includegraphics[width=1.5in]{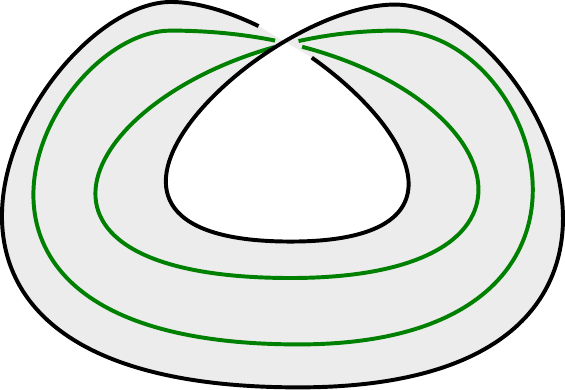}
}
\hspace{1cm}
\subcaptionbox{A relative trisection diagram for $D^2\ttimes\R\bP^2$.
\label{fig:twdiskbundc}}{
\labellist
\endlabellist
\includegraphics[width=1.5in]{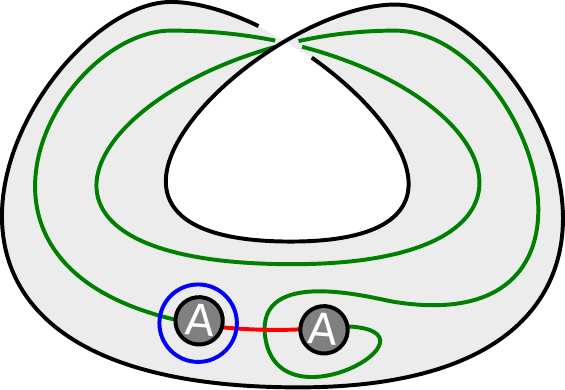}
}

\caption[width=\paperwidth]{From a handle diagram of $D^2\ttimes\R\bP^2$, we obtain a Lefschetz fibration and then a relative trisection diagram.}
\label{fig:d2txrp2}
}
\end{figure}

\begin{example}  \label{ex: twdiskbund} A  handlebody diagram of $D^2 \ttimes  \R\bP^2$, which consists of a single nonorientable $1$-handle and a single $2$-handle  is depicted in Figure~\ref{fig:d2txrp2}(A).  
 We apply the method discussed above to see that $D^2 \ttimes  \R\bP^2$ admits a genus one Lefschetz fibration over $D^2$ with fiber the M\"{o}bius band  and a single vanishing cycle $\gamma$ as shown in Figure \ref{fig:d2txrp2}(B). Since the mapping class group of the M\"{o}bius band is trivial (cf. \cite{e}), the Dehn twist about $\gamma$ is isotopic to the identity and thus  we see that the open book on the boundary $\bdy (D^2 \ttimes  \R\bP^2)$ is trivial, which also implies that  $\bdy (D^2 \ttimes  \R\bP^2) \cong S^2 \ttimes S^1$.
 \end{example}

\begin{figure}[H]{\centering
\subcaptionbox{A handlebody diagram of $D^2\times\R\bP^2$.
\label{fig:twdiskbunda}}{
\labellist
\endlabellist
\includegraphics[width=1.2in]{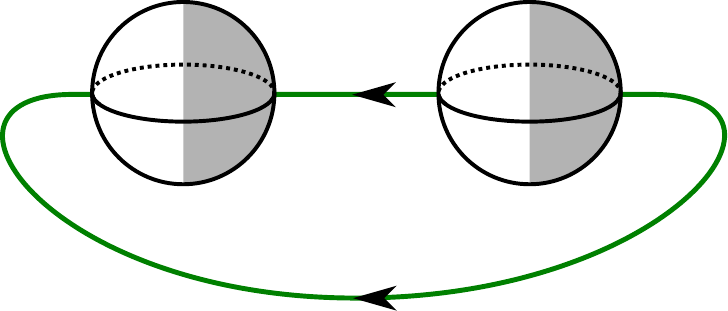}
}
\hspace{.5cm}
\subcaptionbox{We isotope $L$ into a page, but the framing agrees with the page framing.
\label{fig:twdiskbundb}}{
\labellist
\endlabellist
\includegraphics[width=1.2in]{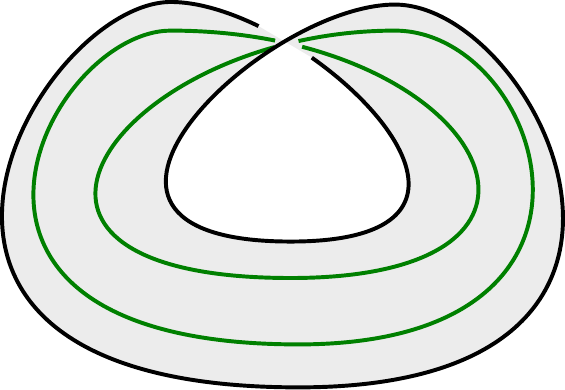}
}
\hspace{.5cm}
\subcaptionbox{$D^2\times\R\bP^2$ admits a genus one Lefschetz fibration over $D^2$ with two vanishing cycles $\gamma_1$ and $\gamma_2$.
\label{fig:twdiskbundc}}{
\labellist
  \pinlabel {\textcolor{darkgreen}{$\gamma_1$}} at 10 0
    \pinlabel {\textcolor{darkgreen}{$\gamma_2$}} at 165 85
\endlabellist
\includegraphics[width=1.2in]{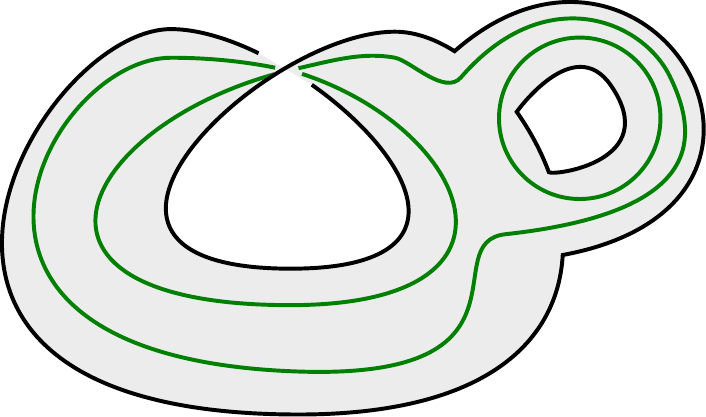}
}
\hspace{.5cm}
\subcaptionbox{A relative trisection diagram for $D^2\times\R\bP^2$.
\label{fig:twdiskbundd}}{
\labellist
\endlabellist
\includegraphics[width=1.2in]{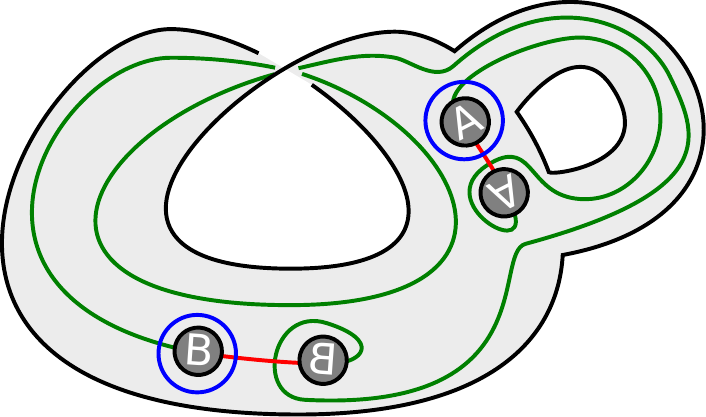}
}

\caption[width=\paperwidth]{From a handle diagram of $D^2\times\R\bP^2$, we obtain a Lefschetz fibration and then a relative trisection diagram.}
\label{fig:prod}
}
\end{figure}

\begin{example} \label{ex: prod}   A handlebody diagram of $D^2 \times  \R\bP^2$, which consists of a single nonorientable $1$-handle and a single $2$-handle  is depicted in Figure \ref{fig:prod}(A). To construct a Lefschetz fibration on $D^2 \times  \R\bP^2$, we need to fix the framing of the $2$-handle by a stabilization. As a result, the fiber will be once stabilized M\"{o}bius band and we will have two vanishing cycle $\gamma_1$ and $\gamma_2$ as shown in Figure \ref{fig:prod}(C). We conclude that $D^2 \times  \R\bP^2$ admits a genus one Lefschetz fibration over $D^2$, whose fiber is $\R\bP^2$ with two holes (i.e., a M\"{o}bius band with one hole) with two vanishing cycles, both of which are boundary parallel. \end{example}

\begin{remark}\label{rem: pi1} Since $\bdy(D^2 \times  \R\bP^2)= S^1 \times \R\bP^2$, as a corollary to  Example~\ref {ex: prod}, we obtain the genus one open book for $S^1 \times \R\bP^2$ whose page is  $\R\bP^2$ with two holes and whose monodromy is the product of Dehn twists about the boundary parallel curves $\gamma_1$ and $\gamma_2$ depicted in Figure~\ref{fig:introfig}, by our construction.  With a natural (but not canonical) choice of local orientations, the Dehn twists will have opposite handedness.  Alternatively, we may consider the abstract open book whose page is  $\R\bP^2$ with two holes and whose monodromy is the product of two boundary Dehn twists,  where we just fix an arbitrary orientation of some annuli neighborhoods of the curves $\gamma_1$ and $\gamma_2$, to specify  the handedness of these Dehn twists. As shown in Figure \ref{fig:dehnsymmetry}, we obtain the same total space regardless of the choice of each handedness. In addition, we can check this fact algebraically: using the method described in \cite[Section 5]{o},  we verify that the fundamental group of the total space of this open book is isomorphic to $\mathbb{Z}\oplus\mathbb{Z}/2\mathbb{Z}$ regardless of the possible choices of handedness for the Dehn twists at hand, as follows:

Choose a basepoint in the boundary of the page $P=\R\bP^2\setminus(\sqcup_2 D^2)$. Let $a$ be a based curve about the generator of $\pi_1(\R\bP^2)$ and let $c_1,c_2$ be based curves about the two boundary components (corresponding to $\gamma_1,\gamma_2$, respectively), so $\pi_1(P)=\langle a,c_1,c_2\mid a^2=c_1c_2\rangle$. See Figure \ref{fig:pi1}.

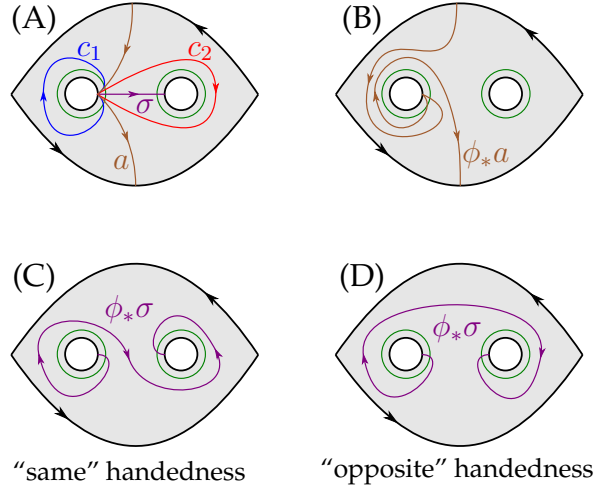
\begin{figure}[H]{\centering
\vspace{.05in}
\labellist
\pinlabel {$\textcolor{blue}{c_1}$} at 50 250
    \pinlabel {$\textcolor{red}{c_2}$} at 120 250
     \pinlabel {$\textcolor{darkbrown}{a}$} at 70 180
     \pinlabel {$\textcolor{darkbrown}{\phi_*a}$} at 300 185
      \pinlabel {$\textcolor{darkpurple}{\sigma}$} at 85 215
        \pinlabel {$\textcolor{darkpurple}{\phi_*\sigma}$} at 75 87
          \pinlabel {\footnotesize``same" handedness} at 75 -15
                \pinlabel {\footnotesize``opposite" handedness} at 282 -15
        \pinlabel {$\textcolor{darkpurple}{\phi_*\sigma}$} at 282 75
         \pinlabel {(A)} at 15 270
            \pinlabel {(B)} at 222 270
               \pinlabel {(C)} at 15 107
                 \pinlabel {(D)} at 222 107
\endlabellist
\includegraphics[width=3in]{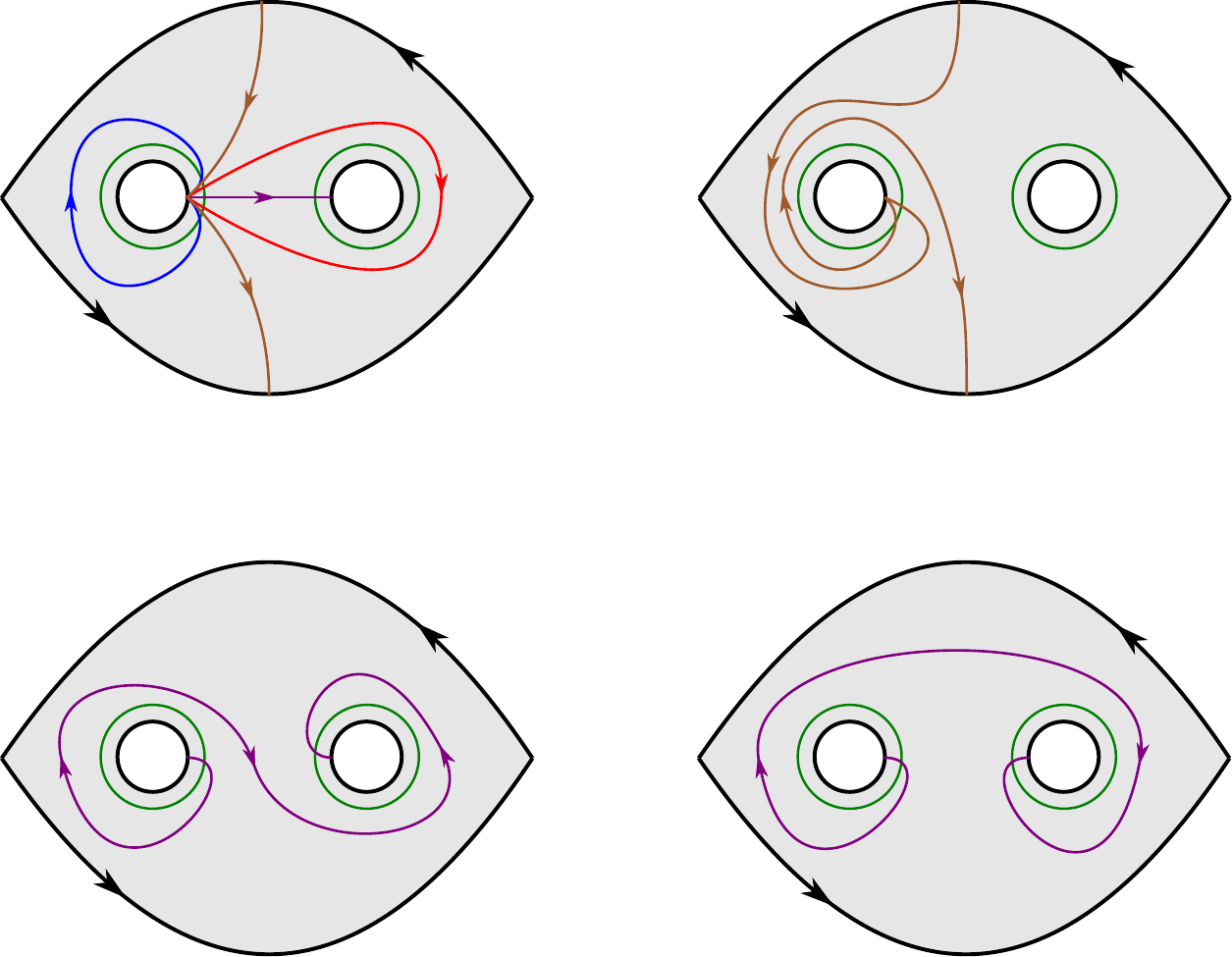}
\vspace{.15in}

\caption{(A) The page $P$ of an open book with monodromy $\phi$ is a copy of $\R\bP^2$ with two disks removed. The map $\phi$ consists of a Dehn twist about each of the green boundary-parallel circles. The twist about the left curve $\gamma_1$ is right-handed with respect to the orientation induced by the ambient page of the figure. (B) We have $\phi_*a=c_1ac_1^{-1}$. (C) If the twist about the right curve $\gamma_2$ is also right-handed, then $(\phi_*\sigma)\sigma^{-1}=c_1c_2^{-1}$. (D) If the twist about $\gamma_2$ is left-handed, then $(\phi_*\sigma)\sigma^{-1}=c_1c_2$. }
\label{fig:pi1}
}
\end{figure}

Now if $\phi$ is the monodromy of an open book with page $P$ and total space $M$, we have \[\pi_1(M)=\langle a,c_1,c_2|a^2=c_1c_2,\phi_*a=a,(\phi_*\sigma)\sigma^{-1}=1\rangle,\] where $\sigma$ is an arc connecting the two boundary components of $P$.

Say that $\phi$ consists of a right-handed (with respect to the local orientation induced by ambient page of Figure \ref{fig:pi1}) Dehn twist about $\gamma_1$ followed by a Dehn twist about $\gamma_2$. Then $\phi_*a=c_1ac_1^{-1}$.

If the Dehn twist about $\gamma_2$ is also right-handed, then $(\phi_*\sigma)\sigma^{-1}=c_1c_2^{-1}$. If the Dehn twist about $\gamma_2$ is left-handed, then $(\phi_*\sigma)\sigma^{-1}=c_1c_2$. In the first case, we obtain
\begin{align*}
\pi_1(M)&=\langle a,c_1,c_2|a^2=c_1c_2,[a,c_1]=1,c_1=c_2\rangle\\
&=\langle a,c_1|a^2=c_1^2,[a,c_1]=1\rangle\\
&\cong\mathbb{Z}\oplus\mathbb{Z}/2\mathbb{Z}.
\end{align*}
In the second case, we obtain
\begin{align*}
\pi_1(M)&=\langle a,c_1,c_2|a^2=c_1c_2,[a,c_1]=1,c_1=c_2^{-1}\rangle\\
&\cong\mathbb{Z}\oplus\mathbb{Z}/2\mathbb{Z}.
\end{align*}

It follows that for any choice of the ``handedness" of a twist about each of $\gamma_1$ and $\gamma_2$, the total space of the resulting open book is $S^1 \times \R\bP^2$.   
\end{remark}

We finish this section with proofs of Corollary~\ref{cor: openb} and Proposition~\ref{prop: minimal}. 

\begin{proof} {\em of Corollary~\ref{cor: openb}.} Let $M$ a closed, connected, nonorientable $3$-manifold. By a theorem of Lickorish \cite{l}, $M$  is obtained from $S^2 \ttimes S^1$ by performing an integral surgery along a link $L$. It follows that $M$ is the boundary of the $4$-manifold $W$ which is obtained by attaching $2$-handles to $D^3 \ttimes S^1$ along $L \subset S^2 \ttimes S^1$. Since $W$ admits a Lefschetz fibration over $D^2$, by Theorem~\ref{thm: lef} the $3$-manifold $\bdy W=M$ acquires an induced open book decomposition, whose monodromy is the product of Dehn twists along the vanishing cycles of the Lefschetz fibration. 
\end{proof}

\begin{proof}{\em of Proposition~\ref{prop: minimal}.} Notice that given an open book of a 3-manifold $M$ with page $P$, one can obtain a Heegaard splitting of $M$ with Heegaard surface equal to the double of $P$. Therefore, if $\#_n S^1\times\R\bP^2$ admits a genus one open book with $k$ binding components, then $\#_n S^1\times\R\bP^2$ admits a Heegaard splitting with Heegaard surface $\Sigma$ the nonorientable surface of genus $2k$.

We make the following observation about Heegaard splittings of nonorientable 3-manifolds: if a Heegaard surface of a nonorientable 3-manifold $N$ is a Klein bottle $\K$, then since $\K$ admits only one nonseparating 2-sided curve (up to isotopy), we must have $N\cong S^2\ttimes S^1$. In particular, every Heegaard surface of $S^1\times\R\bP^2$ has negative (even) Euler characteristic.

Now we again consider our Heegaard surface $\Sigma$ for $\#_n S^1\times\R\bP^2$. By Haken's lemma (see \cite{mn} for some slight discussion about literature in the nonorientable setting), there exist disjoint two-sided curves $C_1,\ldots, C_{n-1}$ embedded in $\Sigma$ so that $\Sigma\setminus \nu(C_1\cup\cdots\cup C_{n-1})$ has $n$ components, each of which closes to a Heegaard surface for a copy of $S^1\times\R\bP^2$. Then \begin{align*}2-2k&=\chi(\Sigma)\\&\le -2n-2(n-1)\\&=2-4n,\end{align*} so we conclude that $k\ge 2n$.
\end{proof}

\section{Trisections of nonorientable Lefschetz fibrations}\label{sec: tri}

In this section, we prove Corollary~\ref{cor: reltri} and draw a trisection diagram for $S^2 \times \R\bP^2$ by doubling a relative trisection diagram of $D^2 \times \R\bP^2$ based on a genus one Lefschetz fibration $D^2  \times  \R\bP^2 \to D^2$ obtained by applying Theorem~\ref{thm: lef}.

\begin{proof} {\em of Corollary~\ref{cor: reltri}. } The construction is essentially the same as in \cite{co}. Consider a Lefschetz fibration over $D^2$ with ordered vanishing cycles $\g_1,\ldots,\g_n$ in a {\em nonorientable} surface $\Sigma$, with total space $X$. Since each $\g_i$ has an annular neighborhood in $\Sigma$, we can perform the wrinkling move on each $\g_i$ as usual. That is, first we replace an annular neighborhood of $\g_1$ with a genus one surface as in Figure \ref{fig:reltri}(B), decorated with one red curve $\alpha_1$, one blue curve $\beta_1$, and one green curve which we also call $\gamma_1$ as an abuse of notation. 

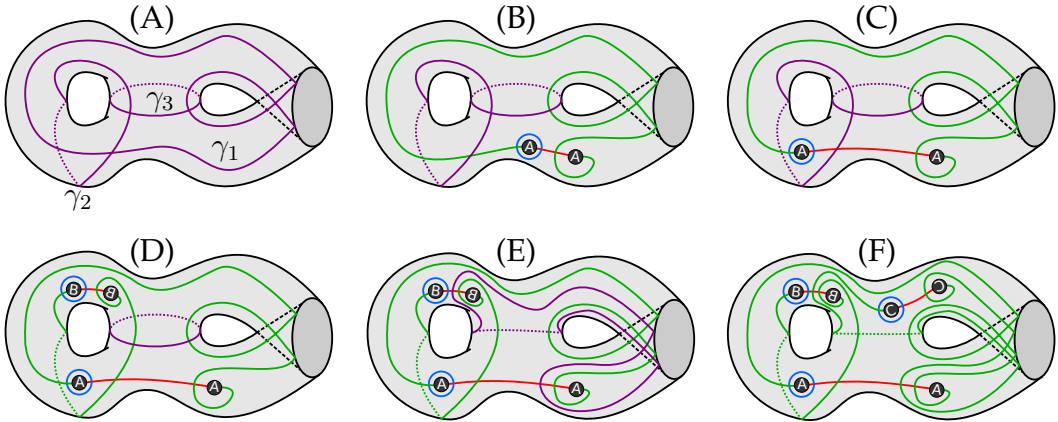
\begin{figure}[H]{\centering
\vspace{.05in}
\labellist
    \pinlabel {$\g_1$} at 135 165
     \pinlabel {$\g_2$} at 45 135
      \pinlabel {$\g_3$} at 95 195
         \pinlabel {(A)} at 90 245
            \pinlabel {(B)} at 312 245
               \pinlabel {(C)} at 532 245
                 \pinlabel {(D)} at 90 103
            \pinlabel {(E)} at 312 103
               \pinlabel {(F)} at 532 103
\endlabellist
\includegraphics[width=5.5in]{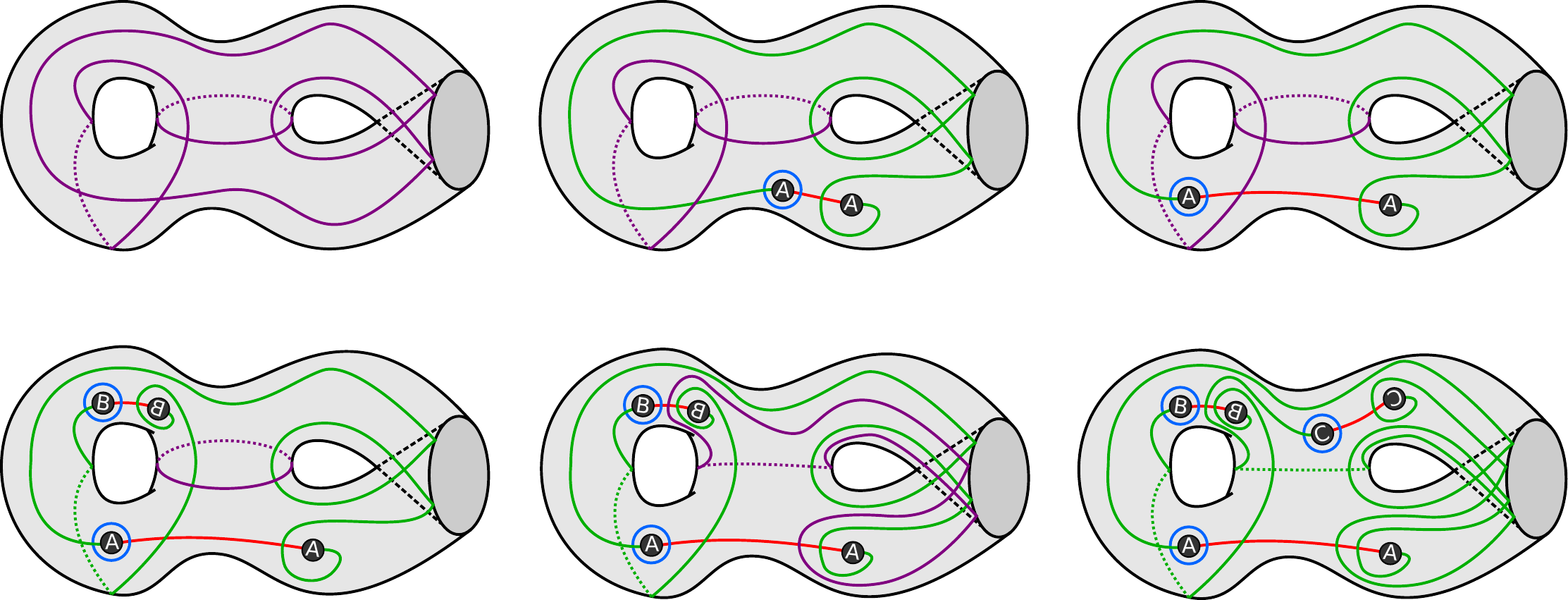}

\caption{In (A), we draw three vanishing cycles $\g_1,\g_2,\g_3$ on a nonorientable surface that together describe a Lefschetz fibration on a nonorientable $4$-manifold $X$. The progression from (A) to (F) illustrates the procedure of Corollary \ref{cor: reltri}. In (F), we obtain a relative trisection diagram of $X$.}
\label{fig:reltri}
}
\end{figure}

Note that in Figure \ref{fig:reltri}(B), $\gamma_1$ appears to twist positively about the tube labels ``A" with respect to a local orientation induced by the page. This sign depends on the ``sign" of the corresponding vanishing cycle: with respect to the same local orientation, it is right-handed near the ``A" tube. If the vanishing cycle were left-handed with respect to this local orientation, then the twist of $\gamma_1$ about ``A" would appear negative.

We then isotope $\g_2$ to avoid $\alpha_1$ and $\beta_1$ and to intersect $\gamma_1$ transversely. Now we slide $\g_2$ over $\beta_1$ as necessary until $\g_2$ is disjoint from $\gamma_1$ (see Figure \ref{fig:reltri}(C)) but now intersects $\alpha_1$. Then we replace an annular neighborhood of $\g_2$ with a genus one surface, again decorated with one red curve $\alpha_2$, a blue curve $\beta_2$, and a green curve $\gamma_2$ (with signs of twisting determined as for $\gamma_1$).  Now isotope $\g_3$ to avoid the red and blue curves and intersect green curves transversely (as in Figure \ref{fig:reltri}(D)), and apply the same procedure to $\g_3$ as to $\g_2$. It is now clear how to iterate this procedure for the rest of the vanishing cycles.

The end result is a surface $\hat{\Sigma}$ obtained from $\Sigma$ by attaching $n$ tubes, and curves $\alpha=\{\alpha_1,\ldots,\alpha_n\},\beta=\{\beta_1,\ldots,\beta_n\},\gamma=\{\gamma_1,\ldots,\gamma_n\}$ on $\hat{\Sigma}$. Then $\mathcal{D}:=(\hat{\Sigma};\alpha,\beta,\gamma)$ is a relative trisection diagram: note that $(\alpha,\beta)$ and $(\beta,\gamma)$ are standard pairs, and since each $\gamma_i$ intersects $\alpha_i$ geometrically once, $(\alpha,\gamma)$ are slide-equivalent to a standard pair. Moreover, $\mathcal{D}$ depicts a $4$-manifold $Y$ obtained from $\Sigma\times D^2$ by attaching $\pm1$-framed $2$-handles along each vanishing cycle in order (see \cite{mn}), so $Y\cong X$. Moreover, via the monodromy algorithm of \cite{cgpc2} (see \cite{mn} about nonorientability), we see immediately that the open book on $\partial X$ induced by $\mathcal{D}$ has page $\Sigma$ with monodromy a Dehn twist on each $\g_i$ in order, as desired. \end{proof}

\begin{example}  Using the genus one Lefschetz fibration $D^2 \ttimes  \R\bP^2 \to D^2$ in Example~\ref{ex: twdiskbund}, we get a  relative trisection diagram for $D^2 \ttimes  \R\bP^2$, by Corollary~\ref{cor: reltri},  as in Figure~\ref{fig:d2txrp2}(C) which coincides with the one depicted  in \cite[Figure 12]{mn}.  Notice that the open book induced on $\bdy(D^2 \ttimes  \R\bP^2)= S^2 \ttimes S^1$ by this relative trisection diagram is trivial, i.e., its page is the M\"{o}bius band and its monodromy is necessarily isotopic to the identity since the mapping class group of the M\"{o}bius band is trivial. Therefore, we can glue this relative trisection diagram for $D^2 \ttimes  \R\bP^2$ with the trivial diagram for $D^3 \ttimes S^1$, to obtain a trisection diagram of $\R\bP^4$, as depicted in \cite[Figure 16]{mn}. \end{example}

\begin{example} \label{ex: s2rp2}  Using the genus one Lefschetz fibration $D^2  \times  \R\bP^2 \to D^2$ in Example~\ref{ex: prod}, we get a relative trisection for  $D^2  \times  \R\bP^2$, by Corollary~\ref{cor: reltri}, as in Figure~\ref{fig:prod}(D).  By doubling  this relative trisection diagram for $D^2  \times  \R\bP^2$,  we get a trisection diagram of $S^2 \times \R\bP^2$ as in Figure~\ref{fig:rp2s2}. 

\begin{figure}[H]{\centering
\vspace{.05in}
\labellist
\endlabellist
\includegraphics[width=3.9in]{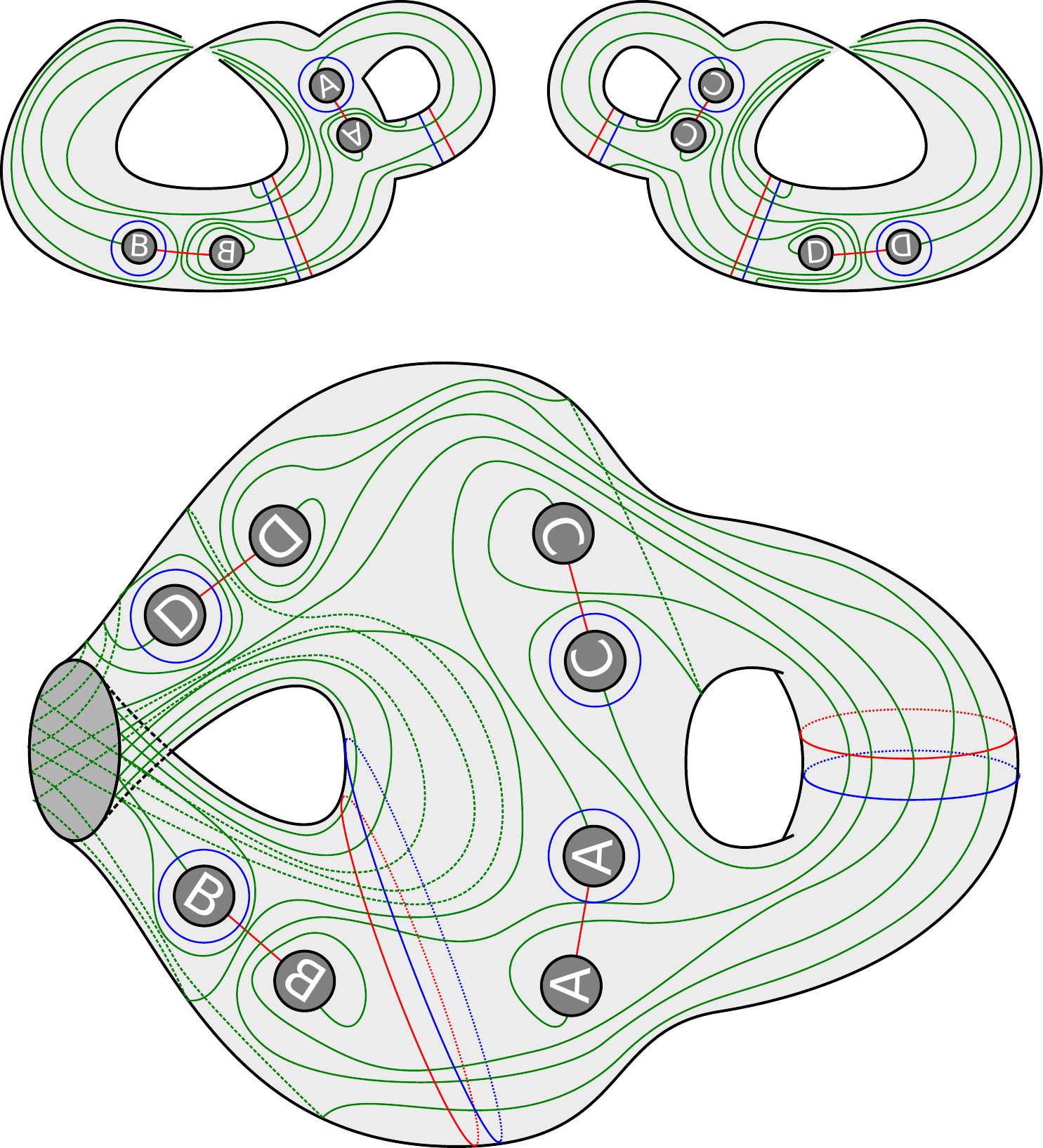}

\caption{We double a relative trisection diagram of $D^2\times\R\bP^2$ to obtain a trisection diagram of $S^2\times\R\bP^2$. {\bf Top:} two copies of the relative trisection of $D^2\times\R\bP^2$ from Figure~\ref{fig:prod}(D). We add extra arcs cutting the $\alpha,\beta,\gamma$ pages into disks in order to perform the gluing operation (see \cite{c}). {\bf Bottom:} gluing the two relative trisection diagrams yields a trisection diagram of $S^2\times\R\bP^2$. Each $\alpha,\beta,\gamma$ arc is doubled to become a new $\alpha,\beta$, or $\gamma$ curve.}
\label{fig:rp2s2}
}
\end{figure}

This example illustrates that these trisections for closed manifolds are generally far from minimal-genus, due in part to the fact that relative trisections have a restrictive boundary condition that often forces high genus. In Figure \ref{fig:marla}, we show a smaller-genus trisection of $S^2\times\R\bP^2$ obtained by applying the procedure of Williams \cite{w} for products of surfaces. (Although Williams only considered orientable surfaces, there is no problem applying her algorithm.)\end{example}

\begin{figure}{\centering
\includegraphics[width=1.5in]{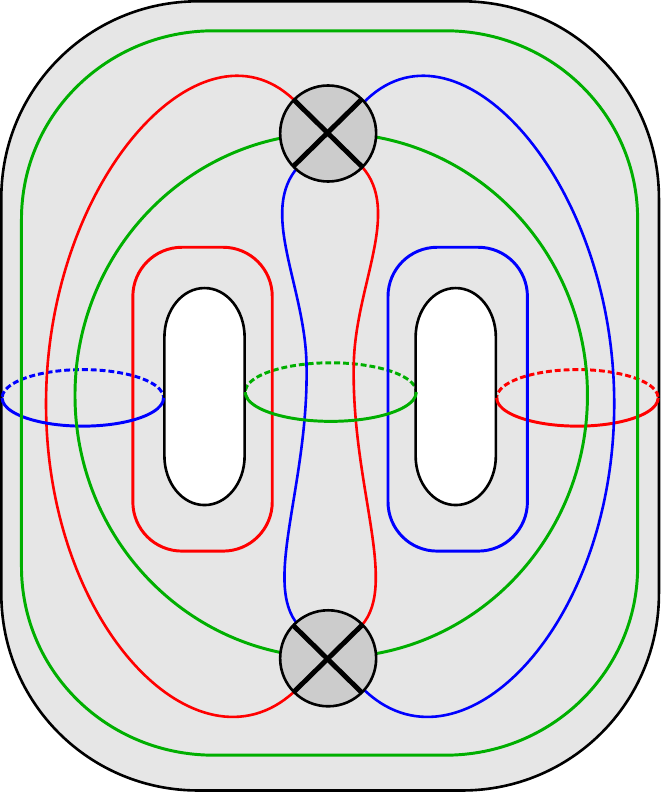}

\caption{A trisection of $S^2\times\R\bP^2$ obtained by following the algorithm of \cite{w}, which produces a trisection surface built by thrice-puncturing two $\R\bP^2$'s and gluing the boundaries together.}
\label{fig:marla}
}
\end{figure}

\begin{remark}\label{rem:yhomeo}
If $X$ is a nonorientable $4$-manifold not including $3$- or $4$-handles with a prescribed open book on its boundary $\del X$, whose  monodromy cannot be factorized into Dehn twists, then Corollary \ref{cor: reltri} cannot be applied to find a relative trisection of $X$ {\em inducing} the given open book. Consider, for example,  the 4-manifold $X$ illustrated in Figure~\ref{fig:yhomeo}(C). It is clear that $X$ admits a Lefschetz fibration over $D^2$ by Theorem \ref{thm: lef}. On the other hand, $\partial X \cong S^1 \times \R\bP^2$ admits an open book with page a Klein bottle with one hole and monodromy the Y-homeomorphism of Lickorish \cite{l}.

\begin{figure}{\centering
\labellist
\pinlabel {$a_1$} at 65 450
\pinlabel {$a_2$} at 175 450
         \pinlabel {(A)} at 106 550
         \pinlabel {$Y(a_1)$} at 310 445
\pinlabel {$Y(a_2)$} at 420 445
            \pinlabel {(B)} at 347 550
               \pinlabel {(C)} at 112 400
                 \pinlabel {(D)} at 349 400
            \pinlabel {(E)} at 222 240
\endlabellist
\includegraphics[width=4in]{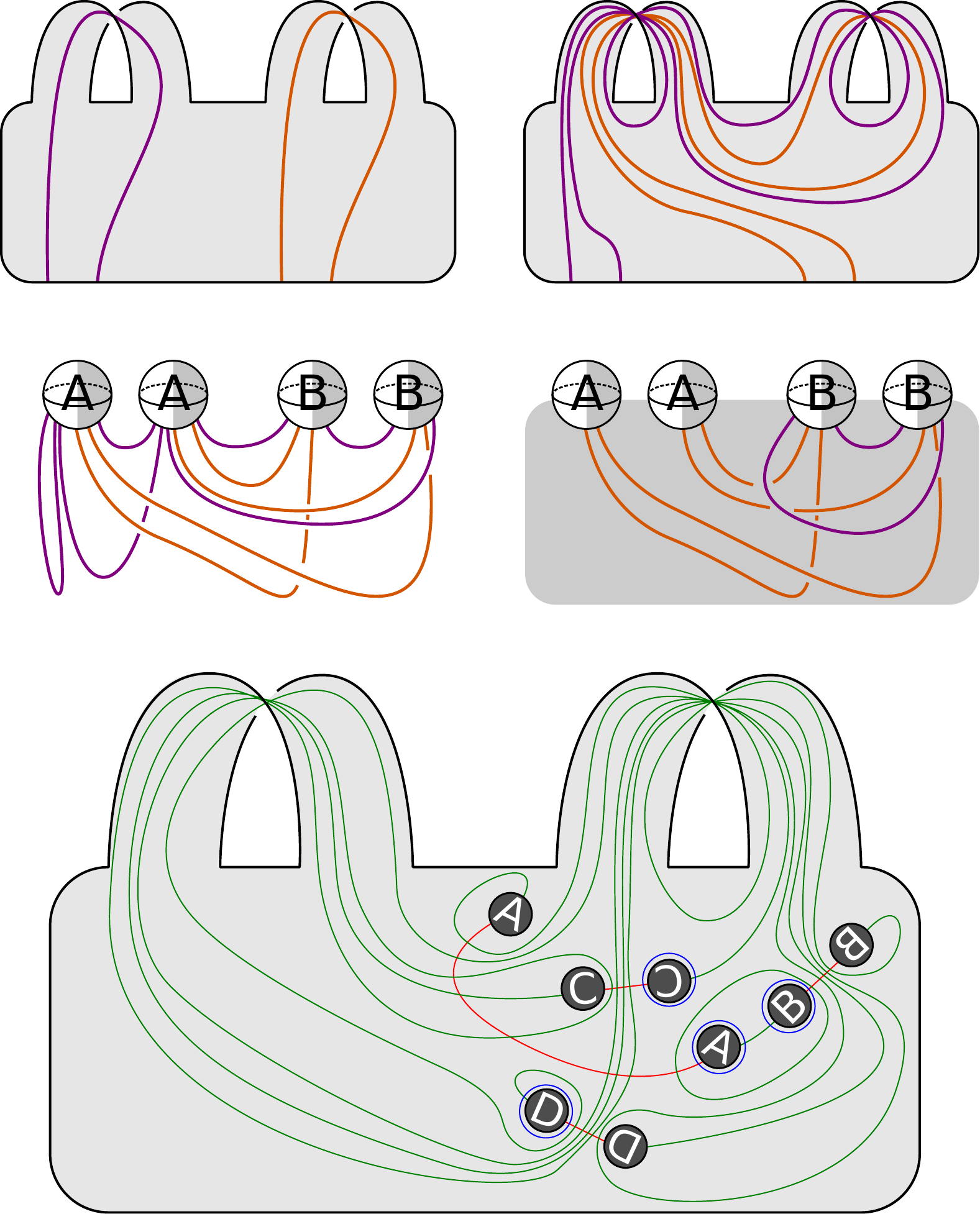}
\vspace{.2in}
\caption{In (A), we show a Klein bottle $\overset{\circ}{\K}$ with one hole and two arcs $a_1,a_2$ whose complement in $\overset{\circ}{\K}$ is a disk. In (B), we show the image of $a_1,a_2$ under the $Y$-homeomorphism; this determines the $Y$-homeomorphism on $\overset{\circ}{\K}$  up to isotopy. In (C), we give a Kirby diagram of a 4-manifold $X$ whose boundary admits an open book $\mathcal{O}$ with page $\overset{\circ}{\K}$ and monodromy the $Y$-homeomorphism. Both $1$-handles are nonorientable and every $2$-handle arc is $0$-framed. In (D), we simplify the diagram slightly and draw a page of the open book. In (E), we obtain a relative trisection of $X$ inducing $\mathcal{O}$ on $\partial X$.}
\label{fig:yhomeo}
}
\end{figure}
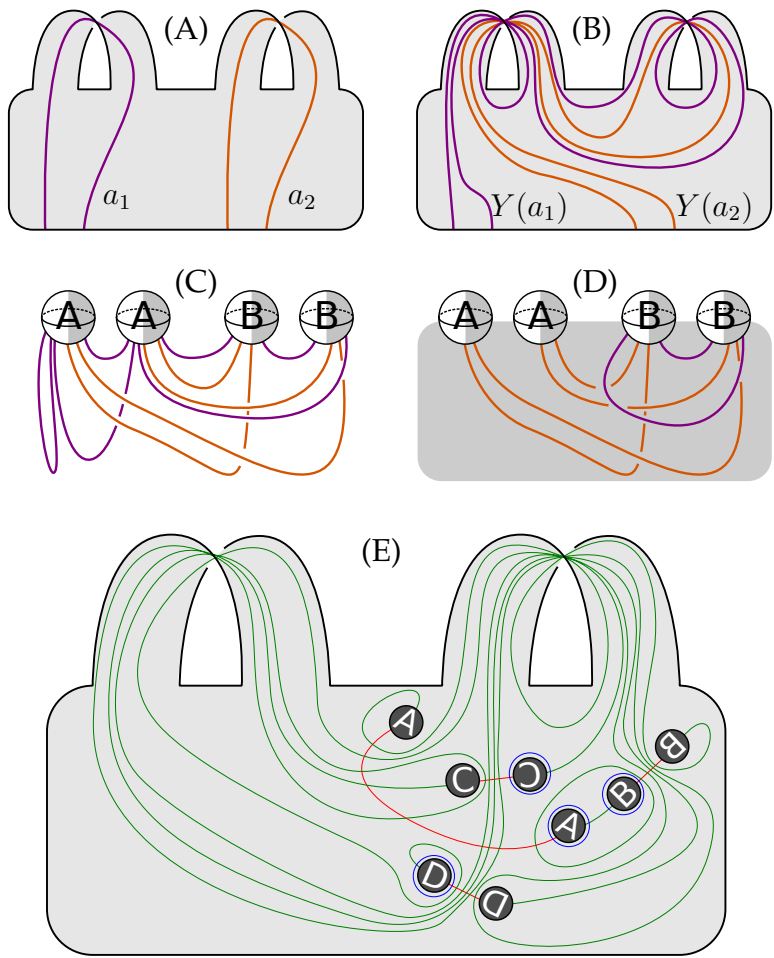

We constructed the handle diagram for $X$ by finding a surgery diagram for the total space of this open book as usual and then replacing all surgery curves with 2-handle attaching circles. This open book extends to a relative trisection  $\mathcal{T}$ of $X$ as illustrated in Figure~\ref{fig:yhomeo}(E) (see \cite{c} for existence, \cite{cgpc} for an algorithm, and \cite{mn} about nonorientability).  As observed by the second author \cite{o}, the open book induced by $\mathcal{T}$ on $\partial X$ is not stably equivalent to an open book whose monodromy is a product of Dehn twists. We conclude that $\mathcal{T}$ is not stably equivalent (via interior and relative trisection boundary stabilization as in \cite{c}) to any relative trisection of $X$ arising from the construction of Corollary \ref{cor: reltri}.
\end{remark}

\begin{remark}
Relative trisections of a 4-manifold $X$ are stably equivalent exactly when the open books they induce on $\partial X$ are related by Hopf plumbing \cite{c}. When $\partial X$ is orientable, then any two relative trisections are related by interior/boundary stabilizations if the open books they induce on $\partial X$ have homologous plane fields \cite{gg} and otherwise by stabilization and one additional move $\partial U$ \cite{cimt}, relying on Harer's  construction \cite{harer} relating any two open books for an orientable 3-manifold by explicit moves. Just like boundary stabilization, the $\partial U$ move has the property that if a relative trisection induces an open book whose monodromy can be written as a product of Dehn twists, then after applying $\partial U$ the new monodromy can still be written as a product of Dehn twists. This means that $\mathcal{T}$ from Remark \ref{rem:yhomeo} cannot be related to any relative trisection of $X$ arising from the construction of Corollary \ref{cor: reltri} by interior/boundary stabilization and $\partial U$ moves, so the classification of relative trisections of a 4-manifold \cite{cimt} does not hold in the nonorientable setting. To classify relative trisections of nonorientable 4-manifolds, it would be necessary to understand relationships between open books on nonorientable 3-manifolds. That is, it would be useful to obtain a move on open books that can transform an open book with monodromy a product of Dehn twists into one whose monodromy cannot be written as a product of Dehn twists, while preserving the total space of the open book.
\end{remark}

A natural question arises as a result of the discussion in Remark~\ref{rem:yhomeo}.

\begin{question}  Does there exist a singular fibration $X \to D^2$ (e.g.  a Lefschetz fibration with ``multiple fibers", a broken Lefschetz fibration or a generic map)  whose regular fiber is the  Klein bottle with one hole and whose monodromy is the $Y$-homeomorphism of Lickorish? 
\end{question} 

\section{Lefschetz fibrations on closed nonorientable $4$-manifolds}\label{sec: closed}

Suppose that $\pi: X \to S^2$ is a Lefschetz fibration on a closed nonorientable $4$-manifold $X$. The topology of $X$ can be described as follows.
First of all, as we pointed out in Remark~\ref{rem: comp}, a regular fiber $F$ is a closed nonorientable surface.   We may also assume that each
singular fiber contains only one critical point, which can be achieved after a small perturbation of $\pi$. Notice that, by the removal of a neighborhood $\nu F$ of a regular fiber $F$, the Lefschetz fibration $\pi: X \to S^2$ canonically determines a Lefschetz fibration $X \setminus \nu F \to D^2$, whose topology is described in Section~\ref{sec: harer}.

Conversely, a Lefschetz fibration (with closed nonorientable fibers) over a disk can be extended to one over a sphere if and only if the monodromy around the disk is trivial. If the genus of the fiber is at least three, this extension is unique by a theorem of Earle and Eells \cite{ee}, which says that  the  group of
diffeomorphisms of any closed nonorientable surface of genus at least three which are homotopic to the identity is contractible.

\begin{remark} \label{rem: pencil} The items listed in \cite[Section 15.2]{os}), each of which is a translation between different types of relations in the mapping class groups of {\em orientable} surfaces and various versions of  {\em orientable} Lefschetz fibrations, hold true verbatim if orientable is replaced with nonorientable. We include two of those items here that will be used below.   Let $\Sigma_{g,k}$ denote the {\em nonorientable} surface of genus $g \geq 1$ with $k \geq 0$ boundary components. A factorization of the identity into a product of Dehn twists along two-sided curves on the closed surface $\Sigma_g$ describes a genus $g$ Lefschetz fibration over $S^2$. Similarly, if $k >0$, a factorization of  $t_{\del_1} \ldots  t_{\del_k}$ (where $\del_i$ is parallel to the $i$th boundary component) into a product of interior Dehn twists along two-sided curves on $\Sigma_{g,k}$ provides a  genus $g$ Lefschetz pencil,  which can be turned into a genus $g$ Lefschetz fibration over $S^2$ equipped with  $k$ disjoint sections, each of square $\pm 1$, by blowing up the base locus. \end{remark}

Notice that when studying Dehn twist factorizations in the mapping class group of a nonorientable surface, one can discard any Dehn twist along a curve that bounds a M\"{o}bius band, since this mapping class is trivial \cite{e}. This is obviously similar to the case of a Dehn twist along a nullhomotopic curve, which is also trivial, whether the surface is orientable or not.  When a factorization contains a Dehn twist along a nullhomotopic curve, this yields a non-relatively minimal Lefschetz fibration, where removing the corresponding singularity---which has no effect on the rest of the fibration---amounts to a surgery that is equivalent to blowing-down a $\C \bP^1$ (possibly with the wrong orientation). If a factorization contains a Dehn twist along a curve that bounds a M\"{o}bius band, one can remove the corresponding singularity in a Lefschetz fibration, which again does not  affect the rest of the fibration. This operation amounts to removing a copy of  $D^2 \ttimes  \R\bP^2$ from the $4$-manifold and gluing in $D^3 \ttimes S^1$. This follows from the fact that the total space of the Lefschetz fibration over $D^2$, whose fiber is the M\"{o}bius band and which has only one boundary-parallel vanishing  cycle is diffeomorphic to $D^2 \ttimes  \R\bP^2$, as we observed in Example~\ref{ex: twdiskbund}, while the trivial (Lefschetz) fibration over $D^2$ whose fiber is the M\"{o}bius band, namely the product of the M\"{o}bius band with $D^2$,  is indeed diffeomorphic to $D^3 \ttimes S^1$. 

Our next  observation is stated as Proposition~\ref{prop: genoneLEF} in the Introduction. \\

\noindent {\bf Proposition~\ref{prop: genoneLEF}.}  {\em  If $X \to B$ is a relatively minimal genus one Lefschetz fibration on a closed nonorientable $4$-manifold $X$ over a closed orientable surface $B$,  then $X$ is an $\R\bP^2$-bundle over $B$.   }

\begin{proof} In a relatively minimal Lefschetz fibration, any vanishing cycle must be homotopically nontrivial in the fiber. But the only homotopically nontrivial simple closed curve in $\R\bP^2$ is one-sided and therefore it can not be a vanishing cycle. We conclude that there are no vanishing cycles in a  relatively minimal {\em genus one} Lefschetz fibration over $S^2$, which proves the desired result. 
\end{proof} 
  
Now we consider relatively minimal genus two Lefschetz fibrations over $S^2$. 

\begin{example}  \label{ex: klein} Lickorish \cite{l} showed that the mapping class group $Map (\K)$ of the Klein bottle $\K$ is isomorphic to $\Z_2 \times \Z_2$, where one copy of $\Z_2$ is generated by the Dehn twist along the only two-sided generic (meaning that it does not bound an embedded disk or a M\"{o}bius band) curve $\a \subset \K$, while the other is generated by the $Y$-homeomorphism. If $t_\a$ denotes the Dehn twist about $\a$, then $t^2_\a$ is isotopic to the identity,  which in turn implies that there exists a {\em genus two} Lefschetz fibration $X(1)  \to S^2$ with only two singular fibers, where $X(1)$ is a closed nonorientable $4$-manifold.  Notice that the extension from a genus two Lefschetz fibration over $D^2$ with monodromy $t_\a^2$, to one over  $S^2$ is not unique (see Remark~\ref{rem: nonunique}) but that is irrelevant for our discussion here. Since there are two singular fibers of $X(1) \to S^2$, it follows that $e(X(1)) = e(K)e(S^2) +2 =2$, where $e$ denotes the Euler characteristic.  Moreover the integral homology $H_1(X(1); \Z)$ is the quotient of $H_1 (K; \Z) = \langle \a, \b \:|\; 2\a=0 \rangle$ by the homology class of  $\a$, and hence it is isomorphic to $\Z$. Next, we define  $X(n)$ inductively as the fiber sum $X(n-1) \#_f X(1)$, inducing a corresponding genus two Lefschetz fibration $X(n) \to S^2$ with monodromy is $t_\a^{2n}$.    It is easy to see that for each positive integer $n$, the $4$-manifold $X(n)$ is nonorientable, $H_1(X(n); \Z)=\Z$,  and $e(X(n))=2n$. Moreover, we can obtain genus two Lefschetz fibrations over an arbitrary  closed orientable surface $B$,  by forming the fiber sum of $X(n)$ with any Klein bottle bundle over $B$. 
\end{example} 

\begin{example} \label{ex: nongeneric}   There exists a  relatively minimal (nonorientable) genus two Lefschetz fibration over $D^2$, with a unique singular fiber whose vanishing cycle is the two-sided essential curve  $\gamma$ that bounds a M\"{o}bius band in $\K$. This Lefschetz fibration can be extended to one over any closed orientable surface $B$, since a Dehn twist about $\gamma$ is isotopic to the identity \cite{e}. \end{example}

\begin{remark} \label{rem: nonunique} The extension from a (nonorientable) genus two Lefschetz fibration over $D^2$ whose monodromy is isotopic to the identity,  to one over $S^2$ is not unique and any two such extensions differ by an element of $\pi_1 (\mbox{Diff}_0 (K))$, which is known to be isomorphic to $\Z$  by \cite{ee}.  Taking this into account, all relatively minimal genus two Lefschetz fibrations over $S^2$ are covered in Examples~\ref{ex: klein} and ~\ref{ex: nongeneric}. \end{remark}

Recall that a logarithmic transformation consists of removing a neighborhood of a torus $T^2$ embedded with trivial normal bundle in an oriented $4$-manifold and gluing back in $T^2 \times  D^2$ by some diffeomorphism of the boundary. The logarithmic transformation has been an essential tool in the classification of minimal elliptic surfaces (see, for instance, \cite[Chapter 8]{gs}.)

\begin{remark}  Similar to generalized log transforms along a torus with trivial normal bundle, it is possible to 
perform surgeries along a Klein bottle with trivial Euler number. In favorable circumstances, such a surgery along the fiber of a nonorientable genus two Lefschetz fibration will preserve the fibration structure. We advise the reader to turn to   \cite{bas} for a discussion of how this surgery is related to the  twisted round $5$-handle attachments and to \cite{n} for its relation to Luttinger surgery on Lagrangian Klein bottles. \end{remark} 

\begin{remark} For any $g \geq 2$,  there exists a  relatively minimal (nonorientable) genus $g$ Lefschetz fibration over any closed orientable surface $B$, with arbitrary number of singular fibers. These can be constructed as in Example~\ref{ex: nongeneric}. \end{remark} 

Next, we turn our attention to finding trisection diagrams of closed nonorientable $4$-manifolds using Lefschetz fibrations. 

\begin{lemma}\label{lem: plumb} Let $V$ denote the $4$-manifold with boundary which is given  by a plumbing of  the disk bundle over $S^2$ with Euler number $\pm 1$ with  the trivial disk bundle over any closed nonorientable surface $F$. Then $V$ admits a Lefschetz fibration over $D^2$  whose regular fiber $\accentset{\circ}{F}$ is obtained by removing a disk from $F$,  which has  only one singular fiber carrying two disjoint vanishing cycles: a homotopically trivial curve and a boundary parallel curve.
\end{lemma}

The statement in Lemma~\ref{lem: plumb} is a very specific case of \cite[Lemma 21]{cgpc2}, where the authors deal with only the plumbings of orientable disk bundles, but the proof goes through verbatim for the  plumbings of nonorientable disk bundles.
Our final result is the following nonorientable analogue of \cite[Theorem 3.7]{co}, which we stated as Theorem~\ref{thm: closedLEF} in the introduction. \\

\noindent {\bf Theorem~\ref{thm: closedLEF}.} {\em Suppose that $X$ is a  closed, nonorientable, connected $4$-manifold  which admits a  Lefschetz fibration over $S^2$  equipped with a section of square $\pm 1$. Then, an explicit  trisection of $X$ can be described  by a  corresponding trisection diagram,  which is determined by the vanishing cycles of the Lefschetz fibration.  }

\begin{proof} Suppose that $X$ is a  closed, nonorientable, connected $4$-manifold and let  $\pi: X \to S^2$  be a Lefschetz fibration, equipped with a sphere section of square $\pm 1$. Let $V$ denote a regular neighborhood of this section union a nonsingular fiber $F$ of $\pi$. Then the $4$-manifold $W$ obtained by removing the interior of $V$ from $X$ admits a  Lefschetz fibration $\pi_W :  W \to D^2$ with fiber  $\accentset{\circ}{F}$, which has  the same set of ordered vanishing cycles as $\pi$.  Note that the total monodromy of the open book (with page  $\accentset{\circ}{F}$) on $\del W $  is isotopic to the boundary Dehn twist $t_{\bdy}$. Applying Corollary~\ref{cor: reltri},  we get a relative trisection on $W$ realizing this open book on $\del W$.

On the other hand, there is a Lefschetz fibration $\pi _V : V \to D^2$ with fiber  $\accentset{\circ}{F}$, as described in Lemma~\ref{lem: plumb}. It follows, by Corollary~\ref{cor: reltri}, that $V$ admits a relative trisection  realizing the open book (with page  $\accentset{\circ}{F}$) on $\del V$ whose monodromy is given by $t_{\bdy}$. Since $\del W= \del V$ by construction, and the open books on $\bdy W$ and $\bdy V$ coincide, we can glue the relative trisection on $W$ with the  relative trisection on $V$, to get a trisection of $X$. In addition, the corresponding relative trisection diagrams can be glued together diagrammatically to obtain a trisection diagram of $X$. It is clear that the final trisection diagram of $X$ will only depend on the vanishing cycles of the initial Lefschetz fibration $\pi: X \to S^2$. \end{proof}

\begin{example}\label{example48}  For any $g >2$,  the boundary Dehn twist $t_\del$ can be expressed as a product of interior Dehn twists along two-sided curves on the surface $\Sigma_{g,1}$. To see this, consider for example the two-holed torus relation $(t_a t_b t_c)^4 = t_d t_e$, and cap-off the boundary component corresponding to $e$ with a nonorientable genus $g-2$ surface with one boundary component. In the resulting nonorientable genus $g>2$ surface with one boundary component, we obtain $t_d = (t_a t_b t_c)^4 t^{-1}_e$. Thus, by Remark~\ref{rem: pencil}, for any $g >2$, there exist a nonorientable genus $g$ Lefschetz fibration over $S^2$ which admits a section of square $\pm 1$. In Figures \ref{fig:section1V}, \ref{fig:section1W}, and \ref{fig:section1} we illustrate the process of obtaining an explicit trisection diagram for the total space of such a Lefschetz fibration for $g=3$, as prescribed by Theorem \ref{thm: closedLEF}. \end{example}

\begin{figure}[H]
\includegraphics[width=5.5in]{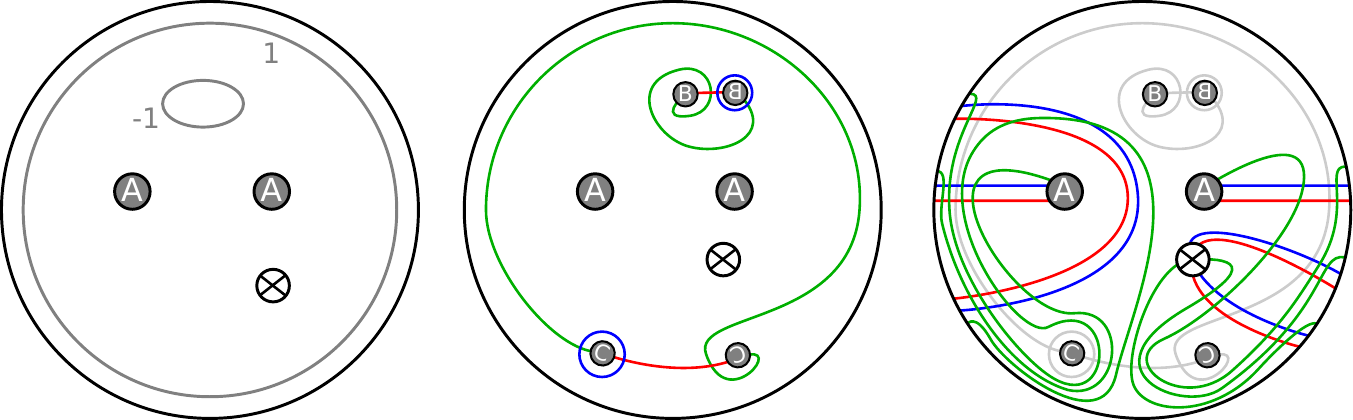}
\caption{Left: the vanishing cycles on the fiber of the Lefschetz fibration $\pi_V:V\to D^2$ as in Theorem \ref{thm: closedLEF}. The manifold $V$ is a neighborhood of of a square one section of $X\to S^2$ along with a regular genus $3$ fiber. There are two vanishing cycles, with signs as indicated. Middle: we obtain a genus $7$, $1$-boundary relative trisection diagram for $V$. Right: cut arcs for the relative trisection diagram. The red and blue arcs are parallel; the green arcs are obtained from these by slides over blue closed curves.}
\label{fig:section1V}
\end{figure}
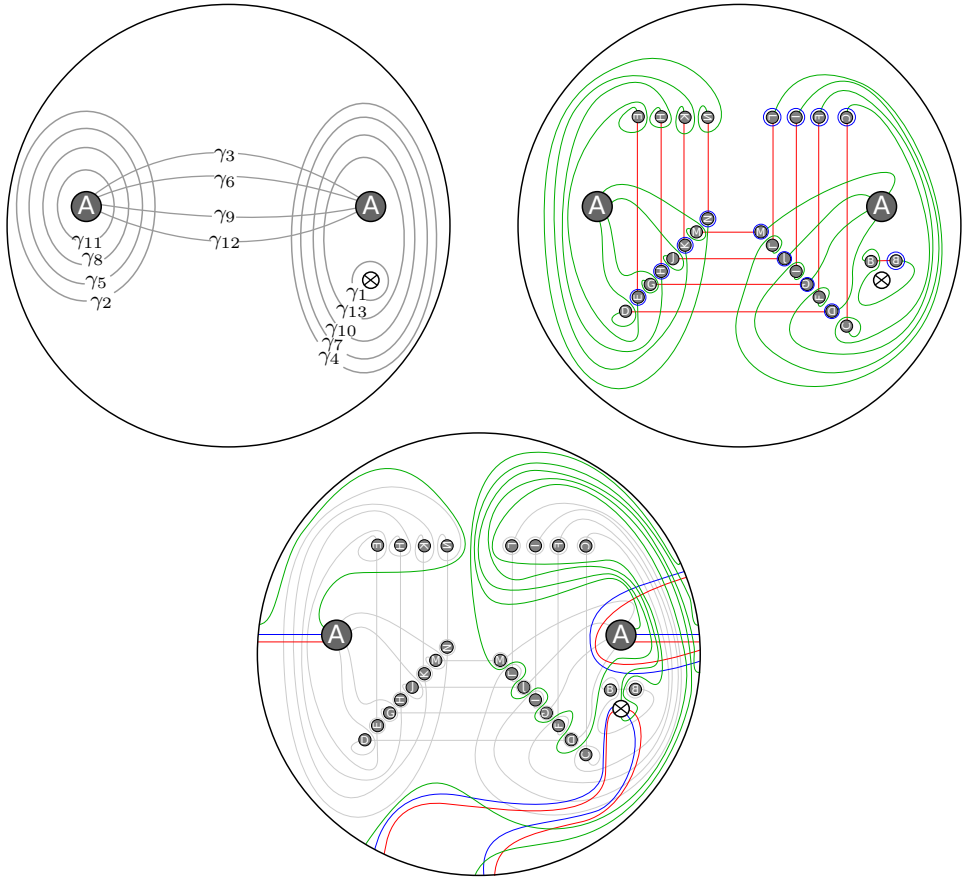
\begin{figure}
\labellist
  \pinlabel {\tiny $\gamma_1$} at 202 335
  \pinlabel {\tiny $\gamma_{13}$} at 199 325
\pinlabel {\tiny $\gamma_{10}$} at 193 314
\pinlabel {\tiny $\gamma_{7}$} at 188 306
\pinlabel {\tiny $\gamma_{4}$} at 186 298
  \pinlabel {\tiny $\gamma_{11}$} at 47 365
\pinlabel {\tiny $\gamma_{8}$} at 50 355
\pinlabel {\tiny $\gamma_{5}$} at 52 342
\pinlabel {\tiny $\gamma_{2}$} at 55 330
  \pinlabel {\tiny $\gamma_{3}$} at 126 415
\pinlabel {\tiny $\gamma_{6}$} at 126 400
\pinlabel {\tiny $\gamma_{9}$} at 126 380
\pinlabel {\tiny $\gamma_{12}$} at 126 365
\endlabellist
\includegraphics[width=5in]{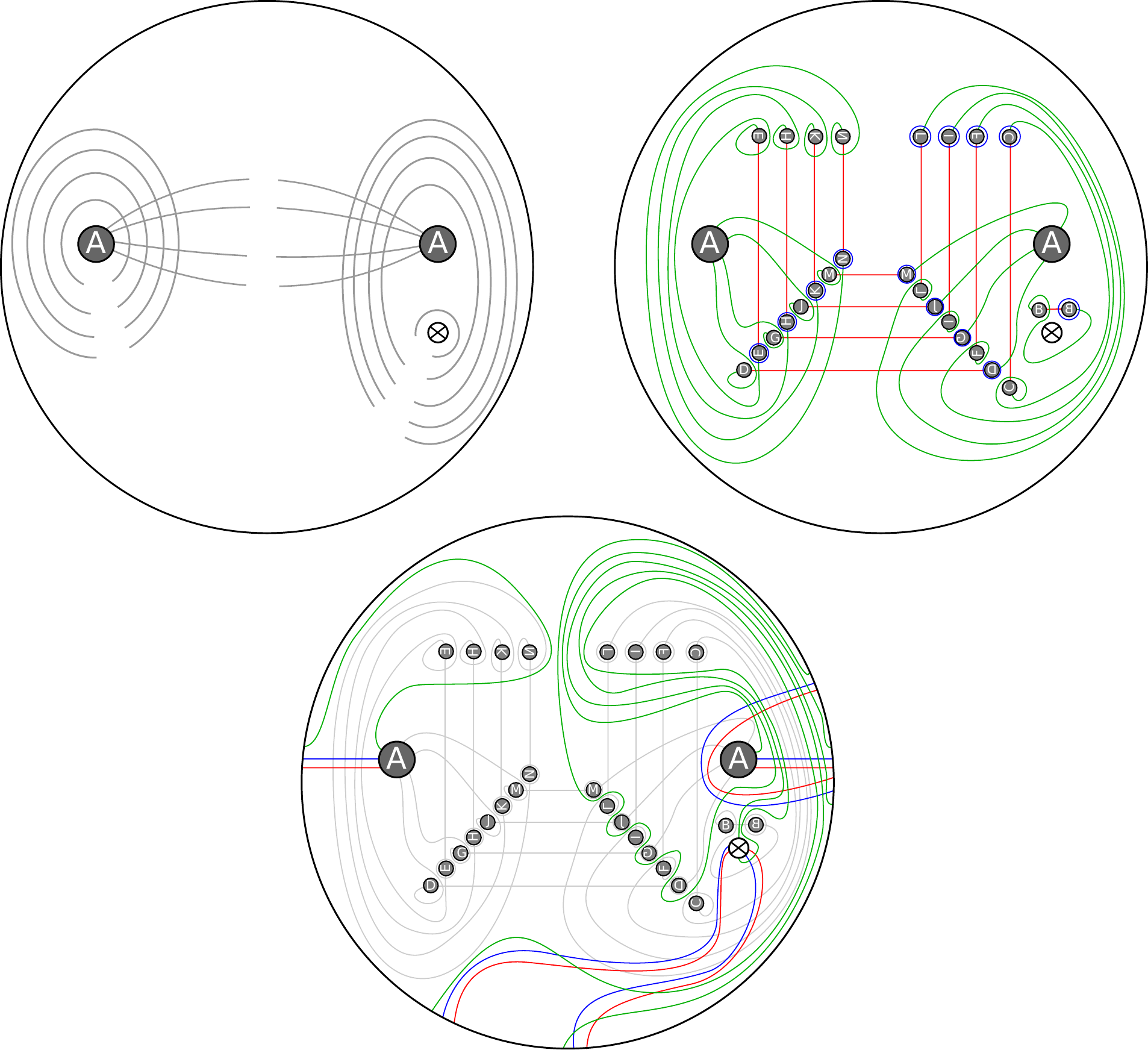}
\caption{Top left: the vanishing cycles on the fiber of the Lefschetz fibration $\pi_W:W\to D^2$ as in Theorem \ref{thm: closedLEF}, where $X$ is described in Example~\ref{example48} with genus $3$ fiber. The vanishing cycles are ordered $\gamma_1,\gamma_2,\ldots,\gamma_{13}$. The vanishing cycle $\gamma_1$ has positive sign while the others have negative sign. Top right: we obtain a genus $29$, $1$-boundary relative trisection diagram for $W$. Bottom: cut arcs for the relative trisection diagram that agree with those in the relative trisection diagram for $V$ from Figure \ref{fig:section1V} after identifying the induced open books on $\partial W$ and $-\partial V$. The red and blue arcs are parallel; the green arcs are obtained from these by slides over blue closed curves.}
\label{fig:section1W}
\end{figure}

In contrast, (nonorientable) genus two Lefschetz fibrations over $S^2$ (for example,  $X(n) \to S^2$ in Example~\ref{ex: klein}) do not admit any sections. This is because there is a unique two-sided generic curve  $\a$ on $\Sigma_{2,1}$ (the Klein bottle with one hole) and for any $m \in \Z$ we have $t_\del \neq t^m_\a$.

\begin{remark} The methods developed in this paper fall short in obtaining an alternate proof, using Lefschetz fibrations, of the fact that every closed nonorientable $4$-manifold admits a trisection analogous to the proof given in \cite[Section 5]{co} for the orientable case, where  contact geometry played an indispensable role.  The obstacle in the nonorientable setting is the nonexistence of an analogue of the Giroux-Goodman stable equivalence theorem \cite[Theorem 1]{gg}.   Nevertheless, it was pointed out to us by the referee that the existence of trisections on nonorientable closed $4$-manifolds can be proven via ``broken" Lefschetz fibrations instead, using the techniques in \cite{bs}. \end{remark}

\begin{figure}[H]
\includegraphics[width=5.5in]{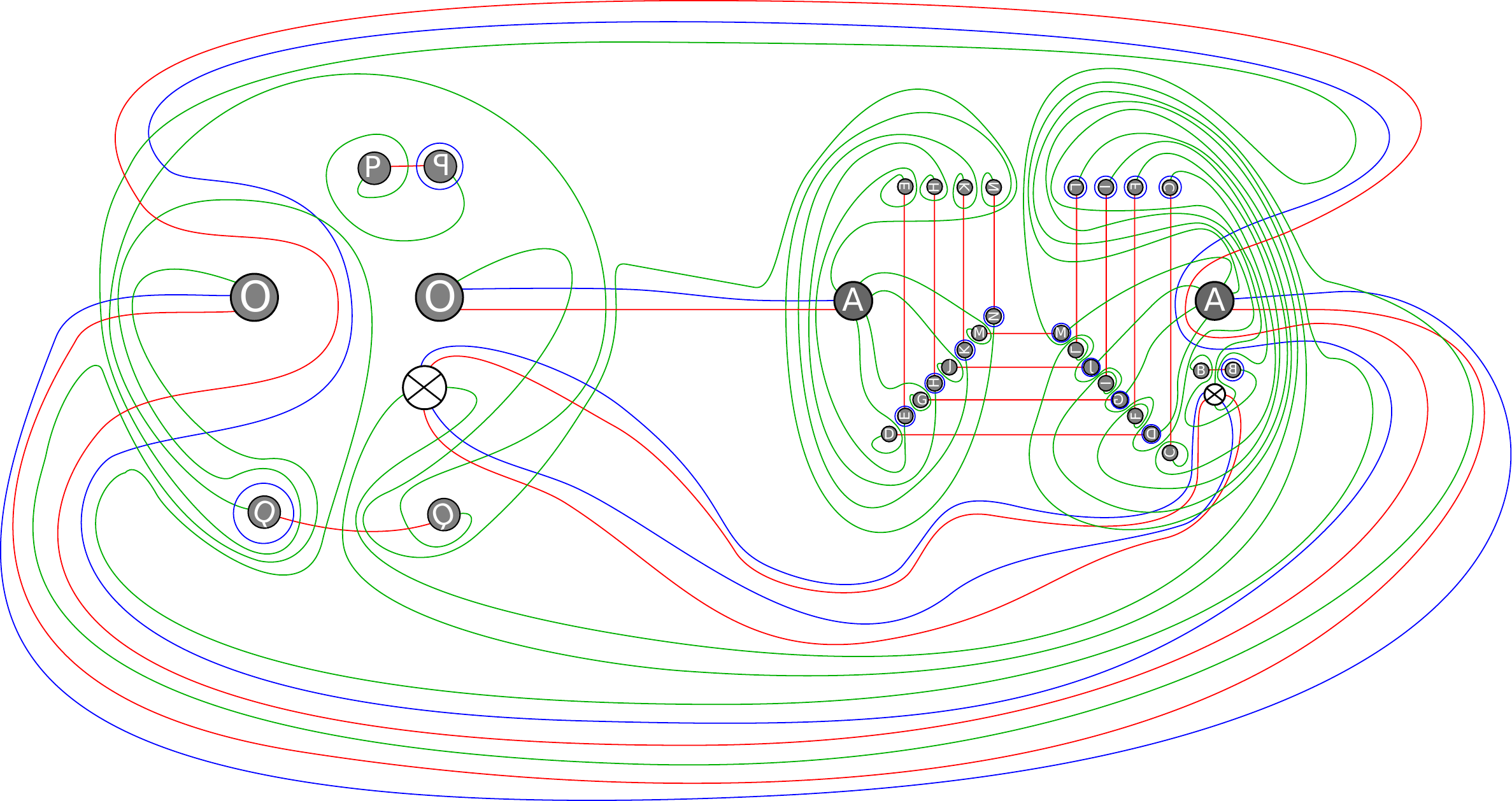}
\caption{We obtain a genus $36$ trisection diagram for $X$ described in Example~\ref{example48} (with genus $3$ fiber) by gluing the relative trisection diagrams for $V$ and $W$, depicted in Figures \ref{fig:section1V} and \ref{fig:section1W}, respectively. The systems of cut arcs glue together to form closed curves in the closed trisection diagram of $X$.}
\label{fig:section1}
\end{figure}

\noindent {\bf {Acknowledgement}}: B.O. would like to thank Mustafa Korkmaz for a useful email correspondence.

\end{document}